\documentclass[12 pt]{article}
\usepackage{amscd, amssymb,amsmath,amsthm}

\usepackage{amsfonts,amssymb,amsmath,eucal,
            amsthm,pstricks,pst-plot,pst-node,eucal}

\usepackage[dvips]{graphics}

\newtheorem{Proposition}{Proposition}

\newtheorem{Lemma}{Lemma}

\newcommand{\proj}{\mathbb{P}}

\newcommand{\rarr}{\rightarrow}
\newcommand{\oh}{{\mathcal{O}}}
\newcommand{\com}{\mathbb{C}}

\newcommand{\eqq}{\stackrel{\sim}{=}}

\def\scup{\mathbin{\text{\scriptsize$\cup$}}}

\newcommand{\bpf}{\noindent {\em Proof.} }
\newcommand{\epf}{\qed \vspace{+10pt}}

\newcommand{\degg}{\text{deg}}
\newcommand{\ds}{\displaystyle}

\begin{document}

\title{Enumerative geometry of Calabi-Yau 4-folds}
\author{A. Klemm and R. Pandharipande}
\date{February 2007}

\maketitle 

\begin{abstract}
Gromov-Witten theory is used
to define an enumerative geometry of curves in Calabi-Yau
4-folds. The main technique is to find exact solutions to 
moving multiple cover integrals. The resulting invariants
are analogous to the BPS counts of Gopakumar and Vafa for 
Calabi-Yau 3-folds. We conjecture the 4-fold invariants 
to be integers and expect a sheaf theoretic explanation.

Several local Calabi-Yau
4-folds are solved exactly.
Compact cases, including 
the 
sextic Calabi-Yau in $\proj^5$, are also studied. 
A complete solution
of the Gromov-Witten theory of the sextic
is conjecturally
obtained by the 
holomorphic anomaly equation.

\end{abstract}

\setcounter{section}{-1}
\section{Introduction}
\subsection{Gromov-Witten theory}
Let $X$ be a nonsingular projective variety over $\com$. 
Let 
$\overline{M}_{g,n}(X,\beta)$ be the moduli space of 
 genus $g$, $n$-pointed stable maps to $X$
representing the class
$\beta\in H_2(X,\mathbb{Z})$.
The Gromov-Witten theory of primary fields{\footnote{We consider
only primary Gromov-Witten theory in the paper.}} concerns the integrals
\begin{equation}\label{ccc}
N_{g,\beta}^X(\gamma_1,\ldots,\gamma_n) =
\int_{[\overline{M}_{g,n}(X,\beta)]^{vir}} \prod_{i=1}^n 
\text{ev}_i^*(\gamma_n),
\end{equation}
where 
$$\text{ev}_i: \overline{M}_{g,n}(X,\beta) \rarr X$$
is the $i^{th}$ evaluation map and
$\gamma_i\in H^*(X,\mathbb{Z})$. The notation
$$N_{g,\beta}^X =
\int_{[\overline{M}_{g}(X,\beta)]^{vir}} 1$$
is used in case there are no insertions.
%When the target $X$ is clear, the superscipt will be dropped.
Since the moduli space 
$\overline{M}_{g,n}(X,\beta)$ is a Deligne-Mumford 
stack, the 
Gromov-Witten invariants \eqref{ccc} are  
$\mathbb{Q}$-valued.

\subsection{Enumerative geometry}
The relationship between  Gromov-Witten theory  and the 
enumerative geometry of curves in $X$
is straightforward in three cases:
\begin{enumerate}
\item[(i)] $X$ is convex (in genus 0),
\item[(ii)] $X$ is a curve,
\item[(iii)] $X$ is $\proj^2$ or $\proj^1 \times \proj^1$.
\end{enumerate}
For (i-iii), the Gromov-Witten theory with primary
insertions equals the classical enumerative
geometry of curves. A discussion of 
convex varieties (i) in genus 0 can
be found in \cite{fp}.
Examples (ii) and (iii) hold
for all genera and recover the \'etale Hurwitz numbers
and the classical Severi degrees respectively.
Case (iii) certainly extends in some form to all
rational surfaces viewed as generic blow-ups. The genus 0
case is treated in \cite{gottp}.

The first nontrivial cases occur for irrational surfaces.
When $X$ is a minimal surface
of general type,
Taubes' results exactly determine the primary
Gromov-Witten invariant for the adjunction genus in the
canonical class,
$$N_{g_X,K_X}^X = (-1)^{\chi(X,\oh_X)},$$
see \cite{t1,t2,t3,t4}.
While much is known about surfaces of general type \cite{park,dm},
surfaces in between are more mysterious. For example, many
questions about the 
relationship of Gromov-Witten theory to the  
enumerative geometry of the $K3$ and Enriques surfaces remain open
\cite{bl,spinning,km,dm}.

The enumerative significance of Gromov-Witten
theory in dimension 3 has been studied since the beginning of the
subject. For Calabi-Yau 3-folds, essentially all
Gromov-Witten invariants, even in genus 0, have large
denominators. The 
Aspinwall-Morrison formula \cite{AM} was conjectured to produce integer
invariants in genus 0. 
A full integrality conjecture for the Gromov-Witten theory of
Calabi-Yau 3-folds 
in
terms of BPS states 
was formulated by Gopakumar and Vafa \cite{GoVafa1,GoVafa2}. Later,
integral invariants 
for all 3-folds were conjectured in \cite{pan,pan2}.  
Various mathematical attempts to
capture the BPS counts in terms of the cohomologies of
associated moduli of sheaves on $X$ were put forward without
a definitive treatment. However, the integral invariants of
\cite{GoVafa1,GoVafa2,pan2} can be conjecturally
interpreted in terms of the sheaf enumeration of Donaldson-Thomas
theory \cite{mnop1,mnop2,t}.

Our main point here is to show the integrality of Gromov-Witten
theory persists in higher dimensions as well. We speculate there
exist universal transformations in every dimension which 
express Gromov-Witten theory in terms of $\mathbb{Z}$-valued invariants.
We conjecture the exact form of the transformation for Calabi-Yau
4-folds. A sheaf theoretic interpretation of the resulting
invariants remains to be found.

\subsection{Calabi-Yau 4-folds} \label{sss}
Let $X$ be a nonsingular, projective, Calabi-Yau 4-fold, and
let $\beta\in H_2(X,\mathbb{Z})$ be a curve class.
 Since 
$$\text{vir dim} \ \overline{M}_{g}(X,\beta) 
= \int_\beta c_1(X) + (\text{dim} \ X -3)(1-g) = 1-g,$$
Gromov-Witten theory vanishes for $g\geq 2$. We need only consider
genus 0 and 1.

We measure the degree of $\beta$ with 
respect to a fixed ample polarization $L$ on $X$,
$$\text{deg}(\beta) = \int_\beta c_1(L).$$
All effective curve classes{\footnote{
Integrality constraints for Gromov-Witten theory always exclude constant
maps. 
The constant contributions
are easily determined in terms of the classical cohomology of $X$. 
For $D\in H^2(X,\mathbb{Z})$,
the genus 1 invariant 
$$N_{1,0}(D)= -\frac{1}{24}\int_X c_3(X) \cup D$$
has denominator bounded by 24.}} 
satisfy $\degg (\beta)>0$.
We abbreviate the latter condition by $\beta>0$.
We are only interested here in Gromov-Witten
invariants for classes satisfying
$\beta>0$.

Integrality in genus 0 is expressed by the following generalization of the
Aspinwall-Morrison formula. 
Invariants $n_{0,\beta}(\gamma_1, \ldots, \gamma_n)$ 
virtually enumerating rational curves of class $\beta$ incident
to cycles dual to the classes $\gamma_i$ are uniquely 
defined by 
\begin{equation}\label{aaa}
\sum_{\beta>0} N_{0,\beta}(\gamma_1,\ldots,\gamma_n)  q^\beta =
\sum_{\beta>0} n_{0,\beta}(\gamma_1,\ldots,\gamma_n) \sum_{d=1}^\infty d^{-3+n} q^{d\beta}.
\end{equation}
A justification for the definition
 via multiple coverings is given in Section \ref{vvv}.

\vspace{+15pt}

\noindent {\bf Conjecture 0}: {\em The invariants $n_{0,\beta}(\gamma_1,\ldots,\gamma_n)$ are integers.}

\vspace{+15pt}

Let $S_1, \ldots, S_s$ be a basis of $H^4(X,\mathbb{Z})$
mod torsion. Let
$$g_{ij}= \int_X S_i \cup S_j $$
be the intersection form, and let
$$\sum_{i,j} g^{ij} [S_i \otimes S_j] \in H^8(X\times X,\mathbb{Z})$$
be the K\"unneth decomposition of the diagonal (mod torsion).

For $\beta_1,\beta_2 \in H_2(X,\mathbb{Z})$, we define 
invariants $m_{\beta_1,\beta_2}$ virtually enumerating rational
curves of class $\beta_1$ meeting rational curves of class $\beta_2$. The meeting
invariants are uniquely determined by the following rules.

\begin{enumerate}
\item[(i)] The invariants are symmetric,
$$m_{\beta_1,\beta_2}=m_{\beta_2,\beta_1}.$$
\item[(ii)] If either $\degg (\beta_1)\leq 0 $ or $\degg(\beta_2)\leq 0$, then $m_{\beta_1,\beta_2}=0$.
\item[(iii)]
If $\beta_1\neq \beta_2,$ then,
$$m_{\beta_1,\beta_2} = 
\sum_{i,j} n_{0,\beta_1}(S_i) \ g^{ij} \ n_{0,\beta_2}(S_j)
 + m_{\beta_1,\beta_2-\beta_1}+
m_{\beta_1-\beta_2,\beta_2}.$$
\item[(iv)]
In case of equality,
$$m_{\beta,\beta} = n_{0,\beta}(c_2(T_X))+
\sum_{i,j} n_{0,\beta}(S_i) \ g^{ij} \ n_{0,\beta}(S_j) -
\sum_{{\beta}_1+{\beta}_2=\beta} m_{\beta_1,\beta_2}.$$
\end{enumerate}
A geometric derivation of the rules (i-iv) is presented in
Section \ref{gen0}.
The conjectural integrality of the invariants 
$n_{0,\beta}(\gamma)$ implies the
integrality of the meeting invariants
$m_{\beta_1,\beta_2}$. 

In genus 1, we need only consider Gromov-Witten
invariants
$N_{1,\beta}$ of $X$
with {\em no} insertions since the virtual dimension is 0.
The invariants $n_{1,\beta}$ virtually
enumerating elliptic curves are uniquely defined by
\begin{eqnarray}\label{bbb}
\sum_{\beta>0} N_{1,\beta} \  q^\beta&  =& \ \ \ \ \ \
\sum_{\beta>0} n_{1,\beta} \ \sum_{d=1}^\infty \frac{\sigma(d)}{d} q^{d\beta} \\ \nonumber & & +\frac{1}{24}
\sum_{\beta>0} n_{0,\beta}(c_2(T_X))\  \log(1-q^\beta)\\ 
\nonumber
&& - \frac{1}{24} \sum_{\beta_1,\beta_2}  m_{\beta_1,\beta_2}\ 
\log(1- q^{\beta_1+\beta_2}).
\end{eqnarray}
The function $\sigma$ is defined by 
$$\sigma(d)= \sum_{i|d} i.$$
The number of automorphism-weighted,
connected, degree $d$, \'etale covers of an elliptic curve is $\sigma(d)/d$.

\vspace{+15pt}

\noindent {\bf Conjecture 1}: {\em The invariants $n_{1,\beta}$ are integers.}

\vspace{+15pt}

\noindent The explicit form of \eqref{bbb} is derived from studying
a particular solvable local Calabi-Yau 4-fold in Section \ref{gen1}.

\subsection{Examples}
The last four Sections of the paper are devoted to the
calculation of basic examples of Calabi-Yau 4-folds.
The two local cases,
$$\oh_{\proj^2}(-1) \oplus \oh_{\proj^2}(-2) 
\rightarrow \proj^2, $$ 
$$
\oh_{\proj^1\times \proj^1}(-1,-1) 
\oplus \oh_{\proj^1\times \proj^1}(-1,-1) \rightarrow \proj^1 \times
\proj^1,$$
are solved in closed form by virtual localization in Section \ref{lc1}.
The local case
$$\oh_{\proj^3}(-4) \rarr \proj^3$$
and the compact Calabi-Yau 4-fold hypersurfaces
$$X_6 \subset \proj^5,$$
$$X_{10}\subset \proj^5(1,1,2,2,2,2), \qquad   X_{2,5}\subset \proj^1\times \proj^4,$$
$$X_{24}\subset \proj^5(1,1,1,1,8,12)$$
are solved by the conjectural holomorphic anomaly equation{\footnote{While
mathematical approaches to the genus 1 invariants in the
compact case are available \cite{gath,zing,dmtop}, the methods are much less
effective than the anomaly equation.}} in Sections
\ref{F1} -\ref{wpp}. The compact cases are much more
interesting than the local toric examples.
In all calculations, the integralities of Conjectures 0 and 1
are verified.

\subsection{Physical interpretation}
Type IIA string compactifications on Calabi-Yau 4-folds give 
rise to massive theories with $(2,2)$ supergravity  
in 2 dimensions. Such theories and their BPS states were 
extensively studied in general~\cite{Cecotti:1992qh,cv} 
and in particular for type IIA on Calabi-Yau 4-folds 
in~\cite{Gates:2000fj,Gukov:1999ya}. The effective
action, worked out in \cite{Gates:2000fj}, 
contains an $\int {\rm d}^2 z R^{(2)}$ term, and the 
topological string at genus 1 calculates a 1-loop 
correction to this term. The latter comes from the 
famous 1-loop term in 10 dimensional type IIA theory that 
was discovered in the context of heterotic type II duality 
in~\cite{Vafa:1995fj} and  gives the following 
contribution to the 10 dimensional effective action 
\begin{equation}
\label{IIAoneloop}
\delta S= -\int d^{10} x \ B \ Y_8(R)\ .
\end{equation}  
Here, $B$ is the $NS-NS$ 2-form of type IIA coupling 
to the string and $Y_8(R)$ is an 8-form constructed 
as a quartic polynomial in the curvature. In 10 dimensions,
the term 
can be directly calculated from the 1-loop amplitude with 4 
gravitons and the antisymmetric $B$-field as external legs. 
If the latter is in the 2 non-compact dimensions,
in the absence of further flux terms, 
the tadpole condition that
$-\frac{\chi(X)}{24}$ vanishes
is obtained. The topological string computes 
the correction to the  $\int {\rm d}^2 z R^{(2)}$ 
term calculated from a loop with 1 external
graviton, 3 internal gravitons, 
and the $B$-field.{\footnote{Work in progress.}} 

As in~\cite{GoVafa1,GoVafa2}, 
the loop integral receives only contributions from BPS states.
The behavior of the topological string amplitude
in the large volume limit 
appears as the zero mass contribution
and supports the claim that the amplitude computes the reduction
of \eqref{IIAoneloop}.
BPS states with a $D2$-brane charge $\beta$ contribute 
$\sim \log(1-q^\beta)$ to the integral. 
The 
integer expansion (\ref{bbb}) can be alternatively written
as
\begin{eqnarray*}
\sum_{\beta>0} N_{1,\beta} \  q^\beta&  =& \ \ \ \ \ \
\sum_{\beta>0} \tilde n_{1,\beta}\ \log(1-q^\beta)
\\ \nonumber & & +\frac{1}{24}
\sum_{\beta>0} n_{0,\beta}(c_2(T_X))\  \log(1-q^\beta)\\ 
\nonumber
&& - \frac{1}{24} \sum_{\beta_1,\beta_2}  m_{\beta_1,\beta_2}\ 
\log(1- q^{\beta_1+\beta_2}).
\end{eqnarray*}
The integrality condition for the invariants
$\tilde{n}_{1,\beta}$ is equivalent to the conjectured
integrality for
 ${n}_{1,\beta}$.
We 
intepret the $\tilde n_{1,\beta}$ as counting BPS states. 
Futhermore, the structure of the  $$ m_{\beta_1,\beta_2}\ 
\log(1- q^{\beta_1+\beta_2})$$
term suggest that the invariants 
$m_{\beta_1,\beta_2}$ count bound states at the threshold of 
BPS states with $D2$-brane charge $\beta_1$ and $\beta_2$
respectively.

\subsection{Outlook}
The meeting invariants make
integrality in genus 1 for Calabi-Yau 4-folds
considerably more subtle than the corresponding
integrality for Calabi-Yau 3-folds. The integrality
transformations in the higher dimensional Calabi-Yau cases should include
all genus 0 meeting configurations. In the non Calabi-Yau
cases, higher genus meeting configurations should occur as well.
Finding the correct coefficients for such a universal
transformation is an interesting problem.

\subsection{Acknowledgments}
The paper began in a conversation with I. Coskun about the classical
enumerative geometry of canonical curves.

We thank M. Aganagic, J. Bryan,
 D. Maulik, S. Theisen, and especially C. Vafa for useful discussions.
A.~K. was partially supported by the  DOE grant DE-FG02-95ER40986.
R.~P. was partially supported by the Packard foundation and
the NSF grant DMS-0500187. The project was started during a visit
of R.~P. to MIT in the fall of 2006.

\section{Genus 0}
\subsection{Aspinwall-Morrison} \label{vvv}
Let $\pi$ and $\iota$ denote 
the universal curve and map over the moduli space,
$$\pi: \mathcal{C} \rarr \overline{M}_{0,0}(\proj^1,d),$$
$$\iota: \mathcal{C} \rarr \proj^1.$$
The Aspinwall-Morrison formula is
\begin{equation*}
\int_{\overline{M}_{0,0}(\proj^1,d)} c_{top}\left(R^1\pi_*\iota^*
(\oh_{\proj^1}(-1)
\oplus \oh_{\proj^1}(-1))\right)= \frac{1}{d^3}.
\end{equation*}
By the divisor equation, we obtain
\begin{equation}\label{jww}
\int_{\overline{M}_{0,n}(\proj^1,d)} c_{top}\left(R^1\pi_*f^*
\left(\oh_{\proj^1}(-1)
\oplus \oh_{\proj^1}(-1)\right)\right) \scup \prod_{i=1}^n \text{ev}_i^*([P])= 
d^{-3+n},
\end{equation}
where $[P]\in H^2(\proj^1,\mathbb{Z})$ is the class
of a point.

Let $X$ be a Calabi-Yau 4-fold, and let
 $$V_1, \ldots, V_n \subset X$$ be cycles imposing a 1-dimensional
incidence constraint for curves. Let
$$C \subset X$$
be a nonsingular rational curve transversely incident to the cycles $V_i$.
If the rational curve has generic normal bundle splitting,
$$N_{X/C} \eqq \oh_{\proj^1}(-1) \oplus \oh_{\proj^1}(-1) 
\oplus \oh_{\proj^1},$$
the contribution of $C$ to the genus 0 Gromov-Witten theory of
$X$ is 
$$\sum_{d =1}^\infty d^{-3+n} q^{d[C]}$$
by \eqref{jww}. The constraints kill the trivial normal direction.

The justification for definition \eqref{aaa} for
the virtually enumerative invariants $n_{0,n}(\gamma_1,\ldots, \gamma_n)$
is complete. Of course, since transversality and genericity
were assumed in the justification, we do not have a proof of
Conjecture 0.

\subsection{Meeting invariants} \label{gen0}
\subsubsection{Rules (i) and (ii)}
Let $X$ be a Calabi-Yau 4-fold.
The meeting invariant $m_{\beta_1,\beta_2}$
virtually enumerates rational 
curves of class $\beta_1$ meeting rational curves of class
$\beta_2$.
Rules (i) and (ii) have clear geometric motivation. In fact, rule (i)
is consequence of rules (ii-iv). Rule (ii) may be viewed
as a boundary condition.

Ultimately, $m_{\beta_1,\beta_2}$ is defined by rules (i-iv).
Rules (iii) and (iv)  are derived by assuming the best possible
behavior for rational curves. 
However, the ideal assumptions are typically false. 
As in Section \ref{vvv},
our derivation can be viewed, rather, as a justification for the definitions.

\subsubsection{Boundary divisor}
For nonzero classes $\beta_1, \beta_2 \in H_2(X,\mathbb{Z})$, let $\bigtriangleup_{\beta_1,\beta_2}$
denote the virtual boundary divisor
$$\bigtriangleup_{\beta_1,\beta_2} \stackrel{\epsilon}
{\rightarrow} \overline{M}_{0,0}(X,\beta_1+\beta_2)$$
corresponding to reducible nodal curves with degree splitting of type
$(\beta_1,\beta_2)$. In the balanced case $\beta_1=\beta_2$, an ordering is taken in
$\bigtriangleup_{\beta_1,\beta_2}$, and $\epsilon$ is of degree 2.

The virtual dimension of $\bigtriangleup_{\beta_1,\beta_2}$ is 0.
Let
$$M_{\beta_1,\beta_2} = \int_{[\bigtriangleup_{\beta_1,\beta_2}]^{vir}} 1 \ \in \mathbb{Q}$$
be the associated Gromov-Witten invariant.
By the splitting axiom of Gromov-Witten theory,
$$M_{\beta_1,\beta_2} = \sum_{i,j} 
N_{0,\beta_1}(S_i)\ g^{ij}\ N_{0,\beta_2}(S_j),$$
following the notation of Section \ref{sss}.
The meeting invariants $m_{\beta_1,\beta_2}$, defined by rules (i-iv),
may be viewed as an integral version of $M_{\beta_1,\beta_2}$.

\subsubsection{Rule (iii)} \label{cxx}

Ideally, the {\em embedded} rational curves in $X$ of class $\beta_i$
occur in complete, nonsingular, 1-dimensional families 
$$F_i\subset \overline{M}_{0,0}(X,\beta_i).$$
Let $\pi_i$ and $\iota_i$ denote 
the universal curve and map over $F_i$,
$$\pi_i:S_i \rarr F_i, $$
$$\iota_i: S_i \rarr X.$$
Since $\beta_1\neq \beta_2$, the families $F_1$ and $F_2$ are distinct.
Ideally, the surfaces $S_i$ are nonsingular and 
the morphisms $\pi_i$ are smooth except for finitely many 1-nodal fibers.
The meeting number $m_{\beta_1,\beta_2}$ is related to the intersection 
\begin{equation}\label{gwwq}
\iota_{1*}(S_1) \cap \iota_{2*}(S_2) \subset X.
\end{equation}
However, the intersection \eqref{gwwq} is not transverse (even ideally).
A fiber of $\pi_1$ may be a component of a reducible fiber 
of $\pi_2$ or vice versa.

The meeting number $m_{\beta_1,\beta_2}$ is defined to count the ideal number of isolated
points of the intersection \eqref{gwwq}. Hence,
$$m_{\beta_1,\beta_2}+ \delta = \int_X \iota_{1*}[S_1]
 \cap \iota_{2*}[S_2] $$
where the correction $\delta$ is determined by the non-transversal intersection
loci.

The number of times a fiber of $\pi_1$ occurs as a component of a reducible fiber
of $\pi_2$ is simply $m_{\beta_1,\beta_2-\beta_1}$.
Similarly, the opposite event occurs $m_{\beta_1-\beta_2,\beta_2}$ times.
The contribution to $\delta$ of each non-transversal is easily determined.
Let
$$C\subset X$$
be a fiber of $\pi_1$ and a component of a reducible fiber of $\pi_2$. Then
$$\delta(C)= \int_C c_1(E)$$
where
$$0 \rarr N_{S_1/C} \oplus N_{S_2/C} \rarr N_{X/C} \rarr E \rarr 0$$
is the normal bundle sequence.
Certainly $N_{S_1/C}$ is trivial and $N_{S_2/C}$ has degree $-1$. By the Calabi-Yau
condition,
$N_{X/C}$ is of degree $-2$. Hence,
$$\delta(C)=-1.$$

Rule (iii) is obtained by expanding the intersection \eqref{gwwq} via the
K\"unneth decomposition of the diagonal. We have
\begin{eqnarray*}
m_{\beta_1,\beta_2} &=& \int_X \iota_*[S_1] \cap \iota_*[S_2] -\delta \\
& = & 
\sum_{i,j} n_{0,\beta_1}(S_i) \ g^{ij} \ n_{0,\beta_2}(S_j)
 + m_{\beta_1,\beta_2-\beta_1}+
m_{\beta_1-\beta_2,\beta_2}.
\end{eqnarray*}

\subsubsection{Rule (iv)}

In case of equality, the meeting number is more subtle. 
%Here, the meeting number should again be thought of
%as the ideal number of isolated points of the self-intersection
%\begin{equation}\label{gwwq}
%\iota(S_\beta) \cap \iota(S_\beta) \subset X.
%\end{equation}
While the surface $S_\beta$ is ideally nonsingular and 
$$\iota: S_\beta \rarr X$$
is ideally an immersion, $\iota$ is not (even ideally) an embedding.
The correct interpretation of $m_{\beta,\beta}$ is twice the number of ideal double points
of $\iota$. The factor of 2 arises from ordering.

The double point formula \cite{fulton}
yields a calculation of $m_{\beta,\beta}$ as a correction
to the self-intersection,
$$m_{\beta,\beta} = \int_X \iota(S_\beta) \cap \iota(S_\beta) 
- \int_{S_\beta} \frac{c(T_X)}{c(T_{S_\beta})}$$
where $c(T_X)$ and $c(T_{S_\beta})$ denote 
the total Chern classes of the respective bundles.
Expanding the correction term (and using the Calabi-Yau condition) we find
$$\int_{S_\beta} \frac{c(T_X)}
{c(T_{S_\beta})} = \int_S c_2(T_X) +c_1(T_{S_\beta})^2
-c_2(T_{S_\beta}).$$
Certainly,
$$n_{0,\beta}(c_2(T_X))= \int_{S_\beta} c_2(T_X).$$
There is a decomposition
$$c_1(T_{S_\beta}) = -\psi + c_1(F_\beta)$$
where $\psi$ is the cotangent line on $S_\beta$ viewed as the 1-pointed space.
Hence
$$\int_{S_\beta} c_1(T_{S_\beta})^2 = \int_{S_\beta} \psi^2  
+ 4 \chi(F_\beta).$$
An elementary geometric argument shows
$$\int_{S_\beta}\psi^2 = 
-\frac{1}{2}\sum_{\beta+\beta_2=\beta} m_{\beta_1,\beta_2}$$
where the right side is number of reducible fibers of $\pi_\beta$.
Since 
$$\int_{S_\beta} c_2(T_{S_\beta}) = 
\chi(S_\beta)=2\chi(F_\beta)+ \frac{1}{2}\sum_{\beta+\beta_2=\beta} m_{\beta_1,\beta_2},$$
only a calculation of the topological Euler characteristic
 $\chi(F_\beta)$ remains.
The formula
$$\chi(F_\beta)=-n_{0,\beta}(c_2(T_X))+     \sum_{\beta+\beta_2=\beta} m_{\beta_1,\beta_2}  $$
is obtained via a Grothendieck-Riemann-Roch calculation applied to the deformation theoretic
characterization of $T_{F_\beta}$,
$$0 \rarr R^0\pi_*(\omega_\pi^\vee) \rarr
R^0\pi_*\iota^* (T_X) \rarr T_{F_\beta} \rarr \oh_{F_\beta}(D) \rarr 0,$$
where $D\subset F_\beta$ is the divisor corresponding to 
nodal fibers of 
$$\pi: S \rarr F_\beta,$$
see \cite{P2} for a similar discussion.

Rule (iv) is obtained by expanding the self-intersection of $\iota_*[S_\beta]$ 
via the
K\"unneth decomposition of the diagonal and putting together the
above surface calculations. We have
\begin{eqnarray*}
m_{\beta,\beta} &=& \int_X \iota_*[S_\beta] \cap \iota_*[S_\beta] 
- \int_S \frac{c(T_X)}{c(T_{S_\beta})}
  \\
& = & 
\sum_{i,j} n_{0,\beta}(S_i) \ g^{ij} \ n_{0,\beta}(S_j)+ n_{0,\beta}(c_2(T_X)) -
\sum_{{\beta}_1+{\beta}_2=\beta} m_{\beta_1,\beta_2}.
\end{eqnarray*}

\section{Genus 1} \label{gen1}
\subsection{Step I}
The equation defining the virtually enumerative invariants
$n_{1,\beta}$,
\begin{eqnarray}
\label{ttt}\sum_{\beta>0} N_{1,\beta} \  q^\beta&  =& \ \ \ \ \ \
\sum_{\beta>0} n_{1,\beta} \ \sum_{d=1}^\infty \frac{\sigma(d)}{d} q^{d\beta} \\ \nonumber & & +\frac{1}{24}
\sum_{\beta>0} n_{0,\beta}(c_2(T_X))\  \log(1-q^\beta)\\ 
\nonumber
&& - \frac{1}{24} \sum_{\beta_1,\beta_2}  m_{\beta_1,\beta_2}\ 
\log(1- q^{\beta_1+\beta_2}),
\end{eqnarray}
is justified in three steps --- one for each term on the right.

The first term is the easiest since the contribution of an
embedded, super-rigid, elliptic curve $E\subset X$ 
of class $\beta$ to
the genus 1 Gromov-Witten theory of $X$ is
$$\sum_{d=1}^\infty \frac{\sigma(d)}{d} q^{d[E]}.$$
See \cite{pan} for a discussion of super-rigidity and
multiple covers of elliptic curves.

\subsection{Step II}\label{twoii}
The second term of \eqref{ttt} is obtained from the contributions
of families of rational curves to the genus 1 Gromov-Witten invariants
of $X$.

Let $F\subset \overline{M}_{0,0}(X,\beta)$
be the ideal nonsingular family of embedded rational curves (as considered
in Section \ref{cxx}). 
We make two further hypotheses.

\vspace{12pt}

\noindent {\bf Condition 1.} {\em
The family $F$ contains {\em no} nodal 
elements.}
\vspace {12pt}

\noindent Hence, the morphism
$$\pi: S \rarr F$$
is a $\proj^1$-bundle.
In fact, Condition 1 is rarely valid, even ideally, and 
will be corrected in Step III.

The contribution of $F$ to $N_{1,d\beta}$ is
expressed as an excess integral over the moduli space of
maps $\overline{M}_{1}(S,d)$ to the $\proj^1$-bundle
representing $d$ times the fiber class.
Let $\widehat{\pi}$ and $\widehat{\iota}$ denote 
the universal curve and map over the moduli space,
$$\widehat{\pi}:C \rarr \overline{M}_{1}(S,d), $$
$$\widehat{\iota}: C \rarr S.$$

\vspace{12pt}

\noindent {\bf Condition 2.} {\em
$R^0\widehat{\pi}_* \widehat{\iota}^*(N_{X/S})$ vanishes.}
\vspace {12pt}

\noindent With the vanishing of Condition 2,
\begin{equation*}
\text{Cont}_F(N_{1,\beta}) = 
\int_{[\overline{M}_{1}(S,d)]^{vir}} c_{top}\left( 
R^1 \widehat{\pi}_* \widehat{\iota}^* N_{X/S}\right).
\end{equation*}

\begin{Lemma} Under the above hypotheses, \label{zzza}
$$\text{\em Cont}_F(N_{1,\beta}) = -\frac{1}{24d} \int_S c_2(T_X).$$
\end{Lemma}

\bpf Consider the relative
moduli space of maps to the fibers of $\widehat{\pi}$,
\begin{equation}\label{vtr}
\overline{M}_{1}(\widehat{\pi},d)\rarr F.
\end{equation} 
We will use the isomorphism
$$\overline{M}_{1}(S,d) \eqq 
\overline{M}_{1}(\widehat{\pi},d).$$
The two virtual classes are easily compared,
\begin{multline}\label{yoyo}
[\overline{M}_{1}(S,d)]^{vir} =\\  -\lambda \cap 
[\overline{M}_{1}(\widehat{\pi},d)]^{vir} +
\chi(F) \cdot [\overline{M}_{1}(\proj^1,d)]^{vir} \in H_*(
\overline{M}_{1}(S,d), \mathbb{Q}).
\end{multline}
On the right, $\lambda$ is the Chern class of the Hodge bundle, 
$\chi(F)$ is the topological
Euler characteristic, and
$\overline{M}_{1}(\proj^1,d)$ a fiber of
\eqref{vtr}.

Consider first the integral
\begin{equation}\label{gqq}
-\int_{[\overline{M}_{1}(\widehat{\pi},d)]^{vir}} \lambda \cdot
c_{top}\left( 
R^1 \widehat{\pi}_* \widehat{\iota}^* N_{X/S}\right).
\end{equation}
Using the basic boundary relation{\footnote{The required
marked point
can be added and removed by the divisor equation.}}
$$\lambda = \frac{1}{12} \bigtriangleup_0 \ \in 
H^2(\overline{M}_{1,1},\mathbb{Q})$$
and the normalization sequence,
we can rewrite \eqref{gqq} as
\begin{equation*}
-\frac{1}{24}\int_{[\overline{M}_{0,2}(\widehat{\pi},d)]^{vir}} 
(\text{ev}_1\times \text{ev}_2)^{*}([\bigtriangleup_{Diag}]) \cdot
\text{ev}_1^*(c_2(N_{X/S})) \cdot
c_{top}\left( 
R^1 \widehat{\pi}_* \widehat{\iota}^* N_{X/S}\right).
\end{equation*}
Finally, using the Aspinwall-Morrison formula, 
$$-\int_{[\overline{M}_{1}(\widehat{\pi},d)]^{vir}} \lambda \cdot
c_{top}\left( 
R^1 \widehat{\pi}_* \widehat{\iota}^* N_{X/S}\right)=
-\frac{d^{-3+2}}{24} \int_S c_2(N_{X/S}).$$

For the second integral, we use the formula from
\cite{GP} for genus 1 contributions,
\begin{multline*}
\chi(F)\int_{[\overline{M}_{1}(\proj^1,d)]^{vir}} 
c_{top}\left( 
R^1 \widehat{\pi}_* \widehat{\iota}^* N_{X/S}\right) = \\
\chi(F)
\int_{\overline{M}_{1}(\proj^1,d)} c_{top}\left(R^1\widehat{\pi}_*
\widehat{\iota}^*
(\oh_{\proj^1}(-1)
\oplus \oh_{\proj^1}(-1))\right)=
 \frac{\chi(F)}{12d}.
\end{multline*}

Summing the first and second integrals, we obtain, by \eqref{yoyo},
\begin{eqnarray*}
\text{Cont}_F(N_{1,\beta})&  = &  \frac{1}{24d}\left( - \int_S c_2(N_{X/S})
+ 2 \chi(F)\right) 
\end{eqnarray*}
Finally, using the Calabi-Yau condition and the geometry
of $\proj^1$-bundles,
\begin{eqnarray*}
c_2(N_{X/S}) & = &c_2(T_X) +c_1(T_S)^2 - c_2(T_S) \\
& = & c_2(T_X)+c_2(T_S) \\
& = & c_2(T_X)+2 \chi(F),
\end{eqnarray*}
concluding the Lemma.
\epf

Modulo the corrections from nodal elements of $F$ to be discussed
in Step III, the derivation of the second term of \eqref{ttt} is complete
since
$$n_{0,\beta}(c_2(T_X)) = \int_S c_2(T_X).$$

\subsection{Step III}
We now relax Condition 1 of Section
\ref{twoii}, but keep Condition 2 in following stronger form.
Let 
\begin{equation} \label{aqqa}
\pi: S \rarr F
\end{equation}
be the ideal family of embedded rational curves of class $\beta$.
Let
$\overline{M}_{1}(S,\widehat{\beta})$ be the
moduli space of maps to $S$
representing a $\pi$-vertical curve class 
$$\widehat{\beta}\in H_2(S,\mathbb{Z}).$$
The
morphism \eqref{aqqa}
is  the blow-up of a $\proj^1$-bundle over finitely many
points corresponding to the 
$$\frac{1}{2}\sum_{\beta_1+\beta_2=\beta} m_{\beta_1,\beta_2}$$
nodal fibers.
Since $\pi$ is not 
 $\proj^1$-bundle, $\widehat{\beta}$ need not
be a multiple of the fiber class. As before
let $\widehat{\pi}$ and $\widehat{\iota}$ denote 
the universal curve and map over the moduli space
$\overline{M}_{1}(S,\widehat{\beta})$.

\vspace{12pt}

\noindent {\bf Condition $2'$.} {\em
$R^0\widehat{\pi}_* \widehat{\iota}^*(N_{X/S})$ vanishes for
every class $\widehat{\beta}$ satisfying 
$$\widehat{\beta}-[\text{\em Fiber}({\pi})]>0.$$}

\vspace {12pt}

The inequality is required in Condition $2'$. 
Ideally, the inequality is violated for connected curves only if 
$\widehat{\beta}$ equals a multiple of a single component of a
reducible nodal fiber of ${\pi}$.
Then,
$$R^0\widehat{\pi}_* \widehat{\iota}^*(N_{X/S})\neq 0.$$
We view $F$ as not contributing at all to the Gromov-Witten
invariants in classes violating the inequality (as these
curves deform away from $F$).

With the vanishing of Condition $2'$,
\begin{equation*}
\text{Cont}_F(N_{1,\widehat{\beta}}) = 
\int_{[\overline{M}_{1}(S,\widehat{\beta})]^{vir}} c_{top}\left( 
R^1 \widehat{\pi}_* \widehat{\iota}^* N_{X/S}\right)
\end{equation*}
for classes $\widehat{\beta}$ satisfying the inequality.

Since the family $F$ may contain nodal elements, Lemma \ref{zzza}
must be modified.
We have
\begin{eqnarray}\label{bob}
\text{Cont}_F\left(\sum_{\beta-[Fiber(\widehat{\pi})] 
>0} N_{1,\beta}\ q^\beta\right)&
= &\ \frac{1}{24} \  n_{0,\beta}(c_2(T_X))\  \log(1-q^\beta)\\ 
\nonumber
&& + \frac{1}{2}
\sum_{d_1=1}^\infty \sum_{d_2=1}^\infty \sum_{\beta_1+\beta_2=\beta} 
{c_{d_1,d_2}}  m_{\beta_1,\beta_2}\ 
                           q^{d_1\beta_1+d_2\beta_2}\, 
\end{eqnarray}
 for universal constants $c_{d_1,d_2}$.
The first term on the right is the uncorrected answer of Lemma \ref{zzza}.
The second term is the correction. The factor of 2 is included
for the double counting induced by the ordering.

The 
universal form of the correction terms follows from the
canonical 
local analytic geometry near the nodal fibers . Let 
$${\pi}^{-1}(p)= E_1 \cup E_2$$
be a nodal fiber.
The local geometry of $S$ near ${\pi}^{-1}(p)$ is
the total space of the node smoothing deformation.
The restriction of $N_{X/S}$  splits
in the form $\oh_S(E_1)\oplus \oh_S(E_2)$.
The universality of the
correction terms then follows.

\begin{Lemma}\label{zxe}
We have
$$\frac{1}{2}
\sum_{d_1=1}^\infty \sum_{d_2=1}^\infty  
{c_{d_1,d_2}} 
                           q_1^{d_1}q_2^{d_2} = -\frac{1}{24} \ \log(1-q_1q_2)
.$$
\end{Lemma} 

Lemma \ref{zxe} concludes Step III and completes the justification
of definition \eqref{ttt} of the invariants
$n_{1,\beta}$.

\subsubsection{Proof of Lemma \ref{zxe}}
By universality, we can prove Lemma \ref{zxe} by considering
any exactly solved geometry that is sufficiently rich to 
yield all the constants $c_{d_1,d_2}$.

The simplest is the following local geometry. Let $S$ be the
blow-up of $\proj^1\times \proj^1$ at the point $(\infty,\infty)$,
$$S= \text{BL}_{(\infty,\infty)}(\proj^1\times \proj^1)
\stackrel{\nu}{\rarr}
\proj^1 \times \proj^1.$$
Let $L_1$ and $L_2$ be line bundles on $S$,
$$L_1= \nu^*\left(\oh_{\proj^1\times\proj^1}(-1,-1)\right), \ \
L_2= \nu^*\left(
\oh_{\proj^1\times \proj^1}(-1,-1)\right)(E),$$
where $E$ is the exceptional divisor.
Let $X$ be the Calabi-Yau total space
$$X= L_1\oplus L_2 \rarr S.$$
Of course, $X$ is non-compact.

The homology $H_2(X,\mathbb{Z})$ is freely spanned
by 
$$H_2(\proj^1\times \proj^1,\mathbb{Z})= \mathbb{Z}[C_1]\oplus
 \mathbb{Z}[C_2]$$ and
$[E]$. 
Let
$$\beta= [C_1] \in H_2(X,\mathbb{Z}).$$ Certainly,
$F_\beta\eqq \proj^1$ and
the associated universal family is
$$\pi:S \rarr \proj^1$$
obtained by composing $\nu$ with the projection onto the second
factor.

The morphism $\pi$ has a unique nodal fiber over $\infty\in F$
which splits as
$$\beta_1=[C_1]-[E],\ \beta_2=[E].$$
Hence, the only nonzero meeting numbers for $X$ are
$$m_{\beta_1,\beta_2}=m_{\beta_2,\beta_1}=1.$$
Condition $2'$ is easily verified for the
family $F$.

\begin{Proposition}\label{fred} We have
$$\text{\em Cont}_F(N_{1,{d_1\beta_1+d_2\beta_2}}) = 
\frac{ \delta_{d_1,d_2}}{12 d_1}$$
for $d_1,d_2>0$.
\end{Proposition}

\bpf
Let $\mathbf{T}^2=\com^* \times \com^*$ act on $\proj^1\times \proj^1$ by
$$(\xi_1,\xi_2)\cdot ([x_0,x_1],[y_0,y_1]) 
= ([\xi_1 x_0,x_1],[\xi_2y_0,y_1])$$
with fixed points
\begin{equation}\label{fff}
(0,0), \ (0,\infty),\ (\infty,0), \ (\infty,\infty).
\end{equation}
The action of $\mathbf{T}^2$ lifts canonically to $S$.
We calculate
\begin{equation}\label{ves}
\text{Cont}_F(N_{1,{d_1\beta_1+d_2\beta_2}}) =
\int_{[\overline{M}_{1}(S,d_1\beta_1+d_2\beta_2)]^{vir}}
c_{top}\left( R^1\widehat{\pi}_* \widehat{\iota}^*( L_1\oplus L_2)\right)
\end{equation}
by $\mathbf{T}^2$-localization.
With the correct $\mathbf{T}^2$-equivariant 
linearizations of $L_1$ and $L_2$, the integral is possible evaluate
explicitly.

Let $s_1$ and $s_2$ denote the weights of the two torus factors of
$\mathbf{T}^2$. 
The tangent weights
of the $\mathbf{T}$-action on
$\proj^1\times \proj^1$ are
$$(-s_1,-s_2), \ (-s_1,s_2), \ (s_1,-s_2), \ (s_1,s_2)$$
at the respective fixed points \eqref{fff}.

\begin{enumerate}
\item[(i)]
Let $\mathbf{T}^2$ act on $\oh_{\proj^1\times \proj^1}(-1,-1)$
with weights
$$s_1+s_2, \ s_1, \ s_2, \ 0$$
at the respective fixed points \eqref{fff}. The choice induces
a canonical $\mathbf{T}^2$-linearization on $L_1$.
\item[(ii)]
Let $\mathbf{T}^2$ act on $\oh_{\proj^1\times \proj^1}(-1,-1)$
with weights
$$s_1, \ s_1-s_2, \ 0, \ -s_2$$
at the fixed points \eqref{fff}. Together with the
canonical linearization on $\oh_S(E)$, the choice induces
a canonical $\mathbf{T}^2$-linearization on $L_2$.
\end{enumerate}

The $\mathbf{T}^2$-localization contributions of the integral
\eqref{ves} over $0\in F$ must first be calculated.
The contribution over $0\in F$ certainly vanishes unless
$d_1=d_2$, 
An unravelling of the formulas shows
\begin{multline*}
\text{Cont}_{0\in F}(N_{1,d\beta_1+d\beta_2})
 = \\
\int_{\overline{M}_{1}(\proj^1,d)} 
\left(\frac{-s_1-\lambda}{s_1}\right) c_{top}\left( -s_1\otimes
\left(R^1\widehat{\pi}_*\widehat{\iota}^*(\oh_{\proj^1}(-1)\oplus
\oh_{\proj^1}(-1))\right)\right).
\end{multline*}
Then, by a straightforward expansion similar to
the proof of Lemma \ref{zzza}, we obtain
the vanishing
$$ \text{Cont}_{0\in F}(N_{1,d\beta_1+d\beta_2})=0.$$

By the vanishing over $0\in F$, the contribution over $\infty\in F$
must be a constant{\footnote{By definition, the contribution
over $\infty$ is a rational function in $s_1$ and $s_2$.}},
$$ \text{Cont}_{\infty\in F}(N_{1,d_1\beta_1+d_2\beta_2})\in \mathbb{Q}.$$
The $\mathbf{T}^2$-action on $S$ has 3 fixed points 
$$p_0,\ p_\infty,\ p'_\infty $$
over
$\infty\in F$. Here, $p_0$ is the fixed point lying
over $(0,\infty)$,
$p_\infty$ is the node of $\pi^{-1}(\infty)$, and
$p'_\infty$ is the remaining fixed point.
With the linearizations
(i) and (ii), $L_1$ 
has weight 0 over $p_\infty,\ p'_\infty$, 
and $L_2$ has weight 0 over $p_0,\ p_\infty$.

Since the $\mathbf{T}^2$-weight of $L_1$ at $p_\infty$ and $p'_\infty$ is 0,
each
node of the fixed map over these produces a
 $\mathbf{T}^2$-trivial factor 
of $R^1\widehat{\pi}_* \widehat{\iota}^*( L_1)$
by the normalization sequence. Each $\mathbf{T}^2$-fixed
component mapping to $\beta_2$  produces a
cancelling $\mathbf{T}^2$-trivial factor of
$R^1\widehat{\pi}_* \widehat{\iota}^*( L_1)$.
Similarly for $L_2$.

The only localization graphs{\footnote{We follow \cite{GP}
for the graphical terminology for the virtual
localization formula}} which survive the $\mathbf{T}^2$-trivial factors
from the 0 weights of $L_1$ and $L_2$
are {\em double combs}.
A double comb is a connected
graph with a single vertex $v_0$ over $p_0$, a single
 vertex $v'_\infty$ over $p'_\infty$, and a single path 
$$v_0 {-} v_\infty {-} v'_\infty.$$
connecting
$p_0$ to $p'_\infty$ through $p_\infty$:

\begin{figure}[!htbp]
  \begin{center}
    \scalebox{0.20}{\includegraphics{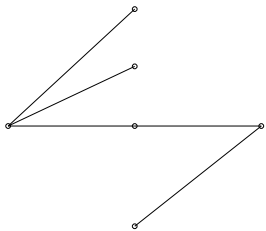}} 
    \caption{A double comb}
    \label{f1}
  \end{center}
\end{figure}

\noindent
Since a double comb has no loops, one of the vertices must
have genus 1. The localization contribution of the double
comb is understood to include all possible genus 1 vertex
assignments.

The final part of the analysis requires taking the nonequivariant
limit 
\begin{equation}\label{nonn}
\text{Lim}_{s_1 \rarr 0} 
\text{Cont}_{\infty\in F}(N_{1,d_1\beta_1+d_2\beta_2}).
\end{equation}
Since the contribution on the right is a constant, no information
is lost.

Nonequivariant limits are often difficult to study, but for
double combs the analysis is simple. 
A factor of $s_1$ in the denominator
of the localization contribution of double comb can only
occur if the two edges of 
the unique path connecting $p_0$ to $p'_\infty$
have equal degrees
$$v_0 \stackrel{e}{-} v_\infty \stackrel{e}{-} v'_\infty.$$
The $s_1$ factor occurs here from the node
smoothing deformation at $v_\infty$.
 Even then,
the $s_1$ factor in the denominator is cancelled
by $s_1$ factors in the numerator if either $v_0$ or $v'_\infty$
has valence greater than $1$.

In case $d_1\neq d_2$, the latter valence condition {\em must}
be satisfied, and the nonequivariant limit \eqref{nonn}
can be taken for each double comb. In fact, each nonequivariant
limit is easily seen to vanish, proving the Proposition in
the unequal case.

If $d_1=d_2$, there is unique double comb which does not
satisfy the valence condition,
\begin{equation}\label{xcff}
v_0 \stackrel{d_1}{-} v_\infty \stackrel{d_2}{-} v'_\infty.
\end{equation}
However, since the nonequivariant limit $\text{Lim}_{s_1\rarr 0}$
exists for {\em all} other double combs, the limit
must exist as well for \eqref{xcff}.
As in the unequal case, the nonequivariant limit
vanishes for all double combs except \eqref{xcff}.
The limit for \eqref{xcff} is explicit calculated
to equal
$$\frac{1}{12d_1}$$
in the equal case.
\epf

To complete the proof of Lemma \ref{zxe}, we expand \eqref{bob}
for the local Calabi-Yau $X$.
Since $c_2(T_X) = 0$,
\begin{eqnarray*}
\text{Cont}_F\left(N_{1,d_1\beta_1+d_2\beta_2}\right)
&=& \frac{1}{2} 
\sum_{d_1=1}^\infty \sum_{d_2=1}^\infty \sum_{\beta_1+\beta_2=\beta} 
{c_{d_1,d_2}}  m_{\beta_1,\beta_2}\ 
                           q^{d_1\beta_1+d_2\beta_2} \\
& = & c_{d_1,d_2}.
\end{eqnarray*}
Hence,
$$c_{d_1,d_2} = \frac{\delta_{d_1.d_2}}{12d}$$
by Lemma \ref{fred}. \qed

The justification for definition \eqref{ttt} of the invariants
$n_{1,\beta}$  is based on ideal geometry. Since the
ideal hypotheses are typically false in algebraic geometry,
Conjectures 0 and 1 are {\em not} proven. In fact, one may be suspicious
of their validity. In the remaining Sections, we will compute many
examples and find the Conjectures to always be valid.
\section{Local examples I} \label{lc1}

\subsection{Solutions}
Proposition \ref{fred} already describes an exactly solved
local Calabi-Yau 4-fold geometry. However, a complete solution is
not given by Proposition \ref{fred} since only special
curve classes of $\text{BL}_{(\infty,\infty)}(\proj^1\times
\proj^1)$ are considered.

The two simplest nontrivial local Calabi-Yau 4-folds are studied here.
The examples may be viewed as the analogues of the
\begin{equation}\label{6yy}
\oh_{\proj^1}(-1) \oplus \oh_{\proj^1}(-1) \rarr \proj^1
\end{equation}
local Calabi-Yau 3-fold. As in \eqref{6yy}, we find
closed form solutions for all curve classes.

\subsection{Local $\proj^2$} 
Let $Y$ be the local Calabi-Yau determined by the total
space of the rank 2 bundle 
$$\oh_{\proj^2}(-1) \oplus \oh_{\proj^2}(-2) 
\rightarrow \proj^2.$$
%$Y$ is the simplest toric Calabi-Yau 4-fold
%with nontrivial Gromov-Witten theory.
Let $P$ denote the point class on $\proj^2$.

\begin{Proposition} We have
$$
N_{0,d}(P) = \frac{(-1)^d}{2d^2} \binom{2d}{d}, $$
\begin{eqnarray*}
\sum_{d>0} N_{1,d} q^{d}  
& = & \frac{1}{12}\log \left(
\sum_{d\geq 0} (-1)^d \binom{2d}{d} q^d
\right)\\
& = & -\frac{1}{24} \log (1+4q). 
\end{eqnarray*}
\end{Proposition}

The Proposition is proven by localization. 
Let $\mathbf{T}^2$ act on $\proj^2$ with fixed points
\begin{equation}\label{nne}
[1,0,0], \ [0,1,0], \ [0,0,1]
\end{equation}
and respective tangent weights
$$(s_1,s_2), \ (-s_2,s_1-s_2), \ (s_2-s_1,-s1).$$
Let $P$ be the equivariant class of the fixed point $[1,0,0]$.
Let $\mathbf{T}^2$ act on $\oh_{\proj^1}(-1)$
with weights
$$0, \ s_2, \ s_1$$
at the respective fixed points \eqref{nne}.
Similarly, let $\mathbf{T}^2$ act on $\oh_{\proj^1}(-2)$ with weights
$$-s_1-s_2, \ s_2-s_1, \ s_1-s_2.$$ 

The above choices kill the localization
contributions to $N_{0,d}(P)$ and
$N_{1,d}$ of all graphs with either a node over $[1,0,0]$ or
an edge connecting $[0,1,0]$ and $[0,0,1]$. The sum over 
remaining comb graphs is not difficult and left to the reader.

The integral invariants $n_{0,d}(P)$ and $n_{1,d}$
can be easily calculated from the Gromov-Witten invariants 
by the defining formulas \eqref{aaa} and \eqref{bbb}.

$$ $$

\begin{tabular}{|c||cccccccccc|}
        \hline
\textbf{}  
&     1  & 2 & 3 & 4 & 5 & 6 & 7 & 8 & 9 & 10  \\
        \hline \hline
$n_{0,d}(P)$  & -1 & 1 & -1 & 2 &-5  &13  & -35 &100 &-300 &  925     \\
$n_{1,d}$      & 0 & 0 & -1 & 2  & -8 & 27 & -90 & 314  &-1140 &4158      \\
       \hline
\end{tabular}

$$ $$

The underlying 
moduli space of maps (with the point condition imposed) for
the invariants $n_{0,1}(P)$ and $n_{0,2}(P)$ are projective
spaces of dimension 1 and 4 respectively. It appears when the
underlying moduli space is $\proj^k$, the invariant is $(-1)^k$
reminiscent of Seiberg-Witten theory for surfaces.

The genus 1 invariants vanish in case there are 
no embedded genus 1 curves. The underlying moduli space for
$n_{1,3}$ is $\proj^9$.

\subsection{Local $\proj^1\times \proj^1$}

Let $Z$ be the local Calabi-Yau determined by the total
space of the rank 2 bundle
$$\oh_{\proj^1\times \proj^1}(-1,-1) 
\oplus \oh_{\proj^1\times \proj^1}(-1,-1) \rightarrow \proj^1 \times
\proj^1.$$
Appropriate localization formulas{\footnote {We now leave the
optimal weight choice for the reader to discover.}}
 for $Y$ in genus 0 and 1
yield 
\begin{multline*}
N_{0,(d_1,d_2)}(P) = \sum_{m\in P(d_1)} \sum_{n\in P(d_2)} 
\frac{(-1)^{d_1+d_2}}{{\mathfrak{z}}(m)
{\mathfrak{z}}(n)} \ \cdot \\
 \int_{\overline{M}_{0,\ell(m)+\ell(n)+1}} 
\frac {1}{ \prod_{i=1}^{\ell(m)} 
(1+m_i \psi_i) \ \prod_{j=1}^{\ell(n)} (1-n_j \psi_{\ell(m)+j})}
\end{multline*}
and
\begin{multline*}
N_{1,(d_1,d_2)}= \sum_{m\in P(d_1)} \sum_{n\in P(d_2)} 
\frac{(-1)^{d_1+d_2}}{{\mathfrak{z}}(m)
{\mathfrak{z}}(n)} \ \cdot \\
 \int_{\overline{M}_{1,\ell(m)+\ell(n)}} 
\frac {1}{ \prod_{i=1}^{\ell(m)} 
(1+m_i \psi_i) \ \prod_{j=1}^{\ell(n)} (1-n_j \psi_{\ell(m)+j})}.
\end{multline*}
Here, $P$ is the point class on $\proj^1\times \proj^1$, and
 $P(d)$ denotes
 the set of partitions of $d$. For $p\in P(d)$, the length is
denoted by $\ell(p)$. The function
$$\mathfrak{z}(p)= |\text{Aut}(p)| \cdot \prod_{i=1}^{\ell(p)} p_i$$
is the usual factor.

By evaluating the above localization sums, we obtain
the following exact solutions.
\begin{Proposition}
We have,
$$N_{0,(d_1,d_2)}(P)=
 \frac{1}{(d_1+d_2)^2}\binom{d_1+d_2}{d_1}^2.$$

\begin{eqnarray*}
\sum_{(d_1,d_2)\neq (0,0)} N_{1,(d_1,d_2)} q_1^{d_1} q_2^{d_2} 
& = & \frac{1}{12}\log \left(
\sum_{d_1\geq 0}\sum_{d_2\geq 0} \binom{d_1+d_2}{d_1}^2 q_1^{d_1} q_2^{d_2}
\right).
\end{eqnarray*}
\end{Proposition}

The integral invariants $n_{0,(d_1,d_2)}(P)$ and $n_{1,(d_1,d_2)}$
can be easily calculated from the Gromov-Witten invariants 
by the defining formulas \eqref{aaa} and \eqref{bbb}.

$$ $$

\begin{tabular}{|c||ccccccc|}
        \hline
\textbf{ $n_{0,(d_1,d_2)}(P)$}  
&  0    &   1  & 2 & 3 & 4 & 5 & 6   \\
        \hline \hline
0   & * & 1 & 0 & 0 & 0 & 0 & 0       \\
1   & 1& 1 & 1 &1  & 1 & 1 & 1        \\
2   & 0 & 1 & 2 & 4 & 6 & 9 & 12          \\
3   & 0 & 1 & 4 & 11 &25 & 49 & 87         \\
4   & 0 & 1 & 6 & 25 & 76 & 196&440        \\
5   & 0 & 1 & 9 & 49 & 196 &635& 1764 \\
6   & 0 & 1 & 12& 87 & 440 &1764 &5926 \\
       \hline
\end{tabular}

$$ $$

\begin{tabular}{|c||ccccccc|}
        \hline
\textbf{ $n_{1,(d_1,d_2)}$}  
&  0    &   1  & 2 & 3 & 4 & 5 & 6   \\
        \hline \hline
0   & * & 0 & 0 & 0 & 0 & 0 & 0       \\
1   & 0 & 0 & 0 &0  & 0 & 0 & 0        \\
2   & 0 & 0 & 1 & 2 & 5 & 8 & 14          \\
3   & 0 & 0 & 2 & 10 &28 & 68 & 144         \\
4   & 0 & 0 & 5 & 28 & 112 & 350 & 922        \\
5   & 0 & 0 & 8 & 68 & 350 &1370& 4426 \\
6   & 0 & 0 & 14& 144 & 922 &4426 & 17220 \\
       \hline
\end{tabular}
$$ $$

For $n_{0,(1,d)}(P)$ the underlying moduli space is $\proj^{2d}$.
The elliptic invariants vanish in classes in which there 
are no embedded elliptic curves.
For $n_{1,(2,2)}$ the moduli space is $\proj^8$.

\section{1-Loop amplitude and Ray-Singer torsion} \label{F1} 
Let $X$ be a nonsingular Calabi-Yau $n$-fold.
The string amplitude which contains information about the 
genus 1 Gromov-Witten theory of $X$ is the twisted 1-loop amplitude 
\begin{equation} 
F_1=\frac{1}{2}\int_{\cal F} \frac{{\rm d}^2 \tau}{ {\rm Im \tau}} {\rm Tr}_{\cal H}\left[ 
(-1)^F F_L F_R Q^{H_L} {\bar Q}^{H_R}\right]\ .
\label{operatordefinition:F1}  
\end{equation}      
Here, the integral is over the fundamental domain ${\cal F}$  of the mapping 
class group of the world-sheet torus with
respect  to the ${\rm SL}(2,\mathbb{Z})$-invariant measure. 
The trace is over the Ramond sector ${\cal H}$ of the twisted 
non-linear $\sigma$-model on $X$.
 The operators $F_L$ and $F_R$ 
are the left and right fermion number operators, $$F=F_L+F_R,$$ 
and $H_L$ and $H_R$ are the left and right moving 
Hamilton operators. The parameter
$\tau$ is the complex modulus of the world-sheet torus, and
$$Q=\exp(2 \pi i \tau ).$$

 The object 
$F_1$ is an index which depends either only on the complexified K\"ahler 
structure moduli of $X$ in the $A$-model or
only 
on the complex structure 
moduli of $\hat X$ in the B-model. The dependence on the 
moduli is via the spectra of $H_L$ and $H_R$. 

We will use the B-model analysis 
to evaluate $F_1$ on the mirror $\hat X$ of $X$.
Predictions 
for the genus 1 invariants of $X$ are then made by the mirror
map. By the world-sheet analysis of \cite{cv,bcov1},
 $F_1$ satisfies the holomorphic 
anomaly equation
\begin{equation} 
\partial_i \bar \partial_{\bar \jmath} F_1=
\frac{1}{2} {\rm Tr}_{\cal H} \left[(-)^F C_i {\bar C}_{\bar \jmath}\right]-
\frac{1}{24}{\rm Tr}_{\cal H}(-)^F\ G_{i\bar \jmath}\ . 
\label{anomaly:F1}  
\end{equation}
Here $G_{i\bar \jmath}$ is the Zamolodchikov metric and 
the derivatives are with respect to the
 $N=2$ moduli. For $N=2$ $\sigma$-models on Calabi-Yau $n$-folds,
 we can specialize to the complex moduli on $\hat X$.
Then, ${\rm Tr}_{\cal H}(-)^F$  
becomes the Euler number $\chi$ of the $\hat X$ and  
$G_{i\bar \jmath}$ becomes 
the Weil-Peterson metric on the complex structure moduli space of $\hat X$. 
The $C_i$ are genus 0, 3-point functions 
in the A-model. 
In the $B$-model on $\hat X$, the $C_i$
 can be calculated from the Picard-Fuchs 
equation for periods of the holomorphic  $(n,0)$ form on $\hat X$. 
The indices $i,\bar \jmath$ run from $1$ to  $h_{n-1,1}(\hat X)$. 

There two methods to integrate the equation (\ref{anomaly:F1}). 
One can use the 
integrability conditions of special geometry for Calabi-Yau $n$-folds 
or, somewhat 
more generally,  the $tt^*$-equations. 
The latter  apply to any 
$N=2$ conformal world-sheet theory. 
If the central charge satisfies
 $$3 c=n \in \mathbb{Z},$$
then the $tt^*$-equations
 imply the special geometry relations for Calabi-Yau $n$-folds. 
The $tt^*$ 
equations are used in \cite{cv,bcov1} to obtain     
\begin{equation}
F_1=\frac{1}{2}\sum_{p,q} (-1)^{p+q} \left(p+q\over 2\right) {\rm Tr}_{p,q} [\log (g)]-\frac{\chi}{24} K +\log |f|^2 \ .
\label{integratedF1}
\end{equation}
The sum here
is over the Ramond-Ramond degenerate lowest energy states labeled by $p,q$ 
which range for the $\sigma$-model case in the left and the right 
moving sector as follows
$$-\frac{n}{2},-\frac{n}{2}+1,\ldots,\frac{n}{2}\, .$$
 By the usual argument \cite{wittenAB}, the states
 are identified in the $A$-model with harmonic 
$$(k,l)=(p+\frac{n}{2},q+\frac{n}{2})$$ 
forms. For 4-folds the $(2,2)$ forms correspond to $(p,q)=(0,0)$ 
and decouple from the sum in (\ref{integratedF1}).
Finally, $g$  is the $tt^*$ metric,
$K$ is the K\"ahler potential for Weil-Peterson metric, and $f$ is the holomorphic 
ambiguity.

% We further note that for the compact
%4-folds discussed below $h_{21}=0$. 
%$g$  denotes the $tt^*$ metric. In the sectors we 
%need it, 

The metric $g$ is related to the Weil-Peterson metric by 
\begin{equation}
G_{i\bar \jmath}=\frac{g_{i\bar\jmath}}{g_{0\bar 0}}=\partial_{i} \bar \partial_{\bar \jmath} K, 
\qquad g_{0\bar 0}=e^{K}\ .
\end{equation} 
The K\"ahler potential $K$ is given by
\begin{equation}
e^{-K}=\int_{\hat X} \Omega \wedge \bar \Omega\ ,
\end{equation} 
where $\Omega$ is the holomorphic $(n,0)$-form on $\hat X$ ---
  $e^{-K}$  can be 
calculated from the periods on $\hat  X$.
    
In summary, specializing to 4-folds{\footnote{All
of our compact examples will satisfy $h_{21}=0$.}} with $h_{21}=0$, we evaluate (\ref{integratedF1}) to
\begin{equation}
F_1=\left(2+h_{11}(X)-\frac{\chi(X)}{24} \right) K -\log \det G+\log |f|^2 \ .
\label{integratedFFF1}
\end{equation}
For 3-folds, in our normalization\footnote{This is, up to the normalization factor ${1\over 2}$, 
the result in~\cite{bcov1}.},
\begin{equation}
F^{(3)}_1=\frac{1}{2}\left(3+h_{11}(X)-\frac{\chi(X)}{12} \right)K -\frac{1}{2}\log \det G+\log |f^{(3)}|^2 \ .
\label{integratedTFF1}
\end{equation}

The Gromov-Witten invariants are extracted in the 
holomorphic limit of (\ref{integratedFFF1}) 
in the large volume of $X$ corresponding 
to the point $P$ of maximal unipotent monodromy of $\hat X$. 
Taking the holomorphic limit is very similar for all dimensions. 
We introduce the flat coordinates near $P$
\begin{equation} 
t_i= \frac{X^i(z)}{X^0(z)}\ ,
\label{mirrormap}
\end{equation}
which are identified with the complexified K\"ahler parameters of $X$. As 
the coordinates
$z$ are the complex structure moduli of $\hat X$,
equation (\ref{mirrormap}) defines
the mirror map between the complex structure on $\hat X$ and the 
complexified K\"ahler structure on $X$. 
The function
$X^0$ is the unique holomorphic period at $P$, which we chose to lie 
at $z_i=0$.
The functions
$$X^i=\frac{1}{2 \pi i}\left( X^0 \log(z)+\text{holomorphic}\right)$$ are the 
$h_{11}(X)=h_{n-1,1}(\hat X)$  single logarithmic periods.
The existence 
of $X^0$ and $X^i$ satisfying the above conditions
 is part of the defining property 
of $P$. Using further the structure of the  periods in an integer 
symplectic basis at $P$, we conclude
\begin{equation}
\begin{array}{rl}
\lim_{\bar t \rightarrow i \infty}K&=-\log(X^0), \\ [3 mm]
\lim_{\bar t \rightarrow i \infty}G_{i \bar \jmath}&= 
\displaystyle{\frac{\partial t_i}{\partial z_j} \delta^j_{\bar \jmath}} \ .
\end{array}
\end{equation} 
After substitution in (\ref{integratedFFF1}), we 
obtain the holomorphic limit of $F_1$ at $P$   
\begin{equation}
{\cal F}_1=\left(\frac{\chi}{24} - h_{11}-2\right)\log(X^0) +
\log \det\left(\frac{\partial z}{\partial t}\right) +\log |f|^2 \ .
\label{holFFF1}
\end{equation}
Here, $\partial z /\partial t$ is the Jacobian of 
the inverse mirror map, and $f(z)$ is the holomorphic ambiguity at genus 1.
The latter is restricted by the
space time modular invariance of $F_1(t,\bar t)$. 
The first two terms in (\ref{integratedFFF1}) can be shown
to be  modular 
invariant. Therefore $f(z)$ must be modular invariant as well.
The modular constraints
together with the large volume behavior, which in physical
terms comes from a zero mode analysis,
\begin{equation} 
{\cal F}_1\rightarrow {(-1)^{n+1}\over 24} \sum_{i=1}^{h_{n-1,1}}
t_i \int_X c_{n-1}(T)\wedge H_i, \qquad t\rightarrow \infty,
\label{limitinfinity}
\end{equation}   
and the expected universal local behavior at other singular limits 
in the complex structure moduli space fix the holomorphic ambiguity 
$f(z)$.

As explained in \cite{bcov}, the the genus 1 free energy  $F_1$ 
is related to the Ray-Singer torsion \cite{raysinger}. The latter
describes aspects of the spectrum of the Laplacians of
$\Delta_{V,q}=\bar \partial_V \bar \partial_V^\dagger+ \bar \partial_V^\dagger \bar \partial_V$ of
a del-bar operator $$\bar \partial_V: \wedge^q \bar T^*\otimes V
\rightarrow \wedge^{q+1} \bar T^*\otimes V$$ coupled to a
holomorphic vector bundle $V$ over $M$. More precisely, with a
regularized determinant over the non-zero mode spectrum  of $\Delta_{V,q}$,
one defines\footnote{\cite{Pestun:2005rp} reviews these facts and relates
the Ray-Singer torsion to Hitchin's generalized 3-form action at one loop.}\cite{raysinger}
\begin{equation} 
I^{RS}(V)=\prod_{q=0}^n \bigl({\det}^{ \prime} \Delta_{V,q}\bigr)^{{q\over 2} (-1)^{q+1}}.
\label{RS-Torsion}
\end{equation}   
The case $V=\wedge^p T^*$ with
$\Delta_{p,q}=\Delta_{\wedge^p T^*,q}$
 leads to the definition of
a family index
\begin{equation} 
F_1={1\over 2}\log \prod_{p=0}^n \prod_{q=0}^n
\bigl( {\det}^{\prime}  \Delta_{pq}\bigr)^{(-1)^{p+q} pq}
\end{equation}   
depending only on  the
complex structure of $\hat X$.

\section{Local examples II} \label{lcn}

We now consider the local Calabi-Yau geometry
$${\cal O}(-n)\rightarrow \mathbb{P}^{n-1}.$$
Since the space is toric,
Batyrev's reflexive cone construction produces the 
mirror geometry: a compact Calabi-Yau $(n-1)$-fold together with
a meromorphic $(n-1,0)$-form  $\lambda$. The latter  can be obtained as a 
reduction of the holomorphic $(n,0)$-form to the Calabi-Yau 
$(n-1)$-fold and has 
a non-vanishing residuum\footnote{One can also consider the elliptic fibration 
over $\mathbb{P}^n$ given by the hypersurface of degree $6n$ in the weighted 
projective space $\mathbb{P}^{n+1}(1^n,2n, 3n)$,
 apply Batyrev's reflexive polyhederal mirror
construction, and take the large fiber limit on both sides.}. The $n$ periods of 
$\lambda$ fulfill the Picard-Fuchs 
equation  
\begin{equation}        
{\cal L} X =\left[\theta^{n-1} - (-1)^n \, n\, z \prod_{k=1}^{n-1} (n\,\theta +k) \right]\theta\, X=0,
\label{mm}
\end{equation}
where $\theta=z \frac{d}{d z}$. The discriminant 
of the Picard-Fuchs equation is $$\tilde \Delta=(1-(-n)^n z).$$

Equation (\ref{mm}) has a constant solution corresponding
 to the residuum of $\lambda$,
 a general property of non-compact Calabi-Yau manifolds.
We normalize constant period to $X^0=1$. 
The system (\ref{mm}) has 3 regular singular points:
\begin{enumerate}
\item[(i)]
the point $z=0$ of maximal unipotent monodromy,
\item[(ii)] the point $\tilde \Delta=0$ corresponding
to a nodal singularity (called the conifold point), 
\item[(iii)]
 the point $z\rightarrow \infty$  (a $\mathbb{Z}_n$ orbifold point). 
\end{enumerate}
Because of (iii), a single 
cover variable $\psi$ is sometimes more convenient. 
It is customary to introduce the latter as
\begin{equation}
z=\frac{(-1)^n}{(n \psi)^n} 
\end{equation}
so the conifold is at $\psi^n=1$. More precisely we define the conifold divisor as  
\begin{equation}
\Delta=(1-\psi^{n})\ .
\end{equation}

Solutions to (\ref{mm}) can be obtained as 
\begin{equation}
X^k=\left. 
\left(\frac{\partial}{2 \pi i \partial \rho}\right)^k  X^0(z,\rho)\right|_{\rho=0}
\label{solutions}
\end{equation}
where we define 
\begin{equation} 
X^0(z,\rho):=\sum_{k=0}^\infty \frac{z^{k+\rho}}{\Gamma(-n (k+\rho)+1) \Gamma(k+\rho+1)^n}  \ .
\end{equation} 

Specializing to $n=4$, we find
 the compact part of the mirror geometry is related to the K3 
given by the quartic in $\mathbb{P}^3$ obtained by setting $n=4$ in 
the above equations.
The meromorphic differential is given by 
$$\lambda= \frac{1}{2 \pi i} \int_{\gamma_0} 
\frac{{\rm d} \Sigma}{p}$$ 
where the contour $\gamma_0$ is around $p=0$ and ${\rm d} \Sigma$ is the
canonical measure on $\mathbb{P}^3$. The  single logarithmic solution is  
\begin{equation} 
X^1=\frac{1}{2 \pi i} \left(\log(z)+ 24\,z + 1260\,z^2 + 123200\,z^3+{\cal O}(z^4)\right)\ .
\label{noncompactmirrormap}
\end{equation} 
We define $q=\exp(2 \pi i t)$ and  with $t=X^1/X^0=X^1$ we obtain by inverting the mirror 
map the series 
\begin{equation} 
z=q - 24\,q^2 - 396\,q^3 - 39104\,q^4 +{\cal O}(q^5)\ .
\end{equation}

The first term of (\ref{holFFF1})  
vanishes in the local case as $X^0=1$.
 The
holomorphic limit of the K\"ahler potential term is trivial. 
We must determine the 
holomorphic ambiguity. As $f$ is a modular invariant,
$f$ can be expressed in terms of $\psi^4$.
As there is a non-degenerate  conformal field theory description at $\psi=0$ given by the 
$\sigma$-model on the orbifold $\mathbb{C}^n/\mathbb{Z}_n$, $F_1$ cannot be singular at this point.
On the other hand the CFT degenerates at $\psi^n=1$ and at $\psi=\infty$, and 
$F_1$ is expected to be logarithmically divergent at the conifold and at the point of 
maximal unipotent monodromy. The former behavior can be argued by comparison with the 
3-fold case while the latter follows directly from (\ref{limitinfinity}) and the leading 
behavior of (\ref{noncompactmirrormap}).
Therefore, we are left with the ansatz $f=x \log(\Delta)$, where $x$ is unknown. We obtain
\begin{equation} 
{\cal F}_1 =\log\left( \frac{\partial \psi}{\partial t}\right) -{1\over 24} \log(\Delta)\ .
\end{equation}   
The first term comes from the holomorphic limit of the Weil-Peterson metric. 
The $x$ coefficient of the last term is matched to the first
 term in the localization calculation 
that can be done in the local case. The leading behavior at boundary
divisors in the moduli space will only depend on the type of the singularities. We expect
therefore that the leading  $-{1\over 24} \log(\Delta)$ behavior will be universal at every
conifold in 4-folds. 
We can check our result for ${\cal F}_1$ against five further terms  
that where calculated applying the localization formulas~\cite{GP} and the Hodge integral formulas~\cite{FP} 
in \cite{mayr}.  For the integer invariants we obtain the results in Table 1. We checked 
integrality of the $n_d$  up to degree $d=100$.

\begin{table}
\begin{centering}
\begin{tabular}{|r|rr|}
\hline
d &g=0 &g=1\\
\hline
1&-20& 0\\ 
2&-820& 0\\ 
3&-68060& 11200\\ 
4&-7486440& 3747900\\ 
5&-965038900& 963762432\\ 
6&-137569841980& 225851278400\\ 
7&-21025364147340& 50819375678400\\ 
8&-3381701440136400& 11209456846594400\\ 
9&-565563880222390140& 2447078892879536000\\ 
10&-97547208266548098900& 531302247998293196352\\ 
11&-17249904137787210605980& 115033243754049262028000\\ 
12&-3113965536138337597215480& 24874518281284024213236000\\  
\hline
\end{tabular}
\caption{Integer invariants  $n_{0,d}$ and $n_{1,d}$   for ${\cal O}(-4)\rightarrow \mathbb{P}^3$ }
\end{centering} 
\end{table}

%-----------------------------------------------------------------

\section{Compact Calabi-Yau 4-folds}\label{wpp}

The holomorphic anomaly equation will now be used verify the 
integrality conjectures for several compact Calabi-Yau 4-folds.
The compact cases have much more interesting geometry 
than the
local models previously considered.

\subsection{The sextic 4-fold} \label{sextic}

We consider first  the sextic $X_6\subset \mathbb{P}^5$. Batyrev's 
mirror construction gives the 1-parameter complex mirror
family  for degree $n$ hypersurfaces  in $\mathbb{P}^{n-1}$ as
\begin{equation} 
p=\sum_{k=1}^n x_k^n- n \psi \prod_{k=1}^n x_k=0\  
\label{polynomial}
\end{equation} 
in\footnote{The orbifold is essentially irrelevant for the B-model period calculation.  It  
only changes the normalization of the periods by a factor $\frac{1}{n^{n-2}}$.} 
$\mathbb{P}^{n-1}/\mathbb{Z}_n^{n-2}$.  The holomorphic  $(n,0)$-form  can be 
written as $$\Omega=\frac{1}{2 \pi i} \int_{\gamma_0} 
\frac{\psi {\rm d} \Sigma}{p}$$
 where the contour $\gamma_0$ is around $p=0$. 
We obtain the Picard-Fuchs operator for the period integrals $\int_{\Gamma^{(4)}} \Omega$ in 
this family parameterized by the variable $z=(-n \psi)^{-n}$ as 
\begin{equation}
{\cal L}= {\theta}^{n-1} - n\,z \prod_{k=1}^{n-1} \left( n\,\theta +k\right)\ .
\end{equation}
The discriminant is $\tilde \Delta=1-n^n z$, and $z=0$ is the 
point of maximal unipotent monodromy. Solutions as in  
(\ref{solutions}) are obtained from 
\begin{equation} 
X^0(z,\rho):=\sum_{k=0}^\infty \frac{\Gamma(n(k+\rho)+1)}{\Gamma(k+\rho+1)^n} z^{k+\rho}\ .
\end{equation} 
Here, there is a non-trivial holomorphic solution
\begin{equation} 
X^0=\sum_{k=0}^\infty \frac{(n k)!}{(k!)^n} z^k\ . 
\end{equation} 
 For the sextic,  the first few terms of the 
inverse mirror map are
\begin{equation} 
 z=q - 6264\,q^2 - 8627796\,q^3 - 237290958144\,q^4 + {\cal O}(q^5)\ . 
\end{equation}
 
The nonvanishing Hodge numbers of $X_6$ up to the symmetries of the Hodge diamond are
$$h_{00}=h_{4,0}=1, \ \  h_{11}=1, \ \ h_{31}=426, \ \ h_{22}=1752.$$ 
Further one has
\begin{equation} 
\chi=\int_X c_4=2610,\quad c_2=15 H^2, \quad c_3=-70 H^3, \quad \int_X H^4=6\ 
\end{equation} 
where $H$ is the hyperplane class.
 
The holomorphic ambiguity can be fixed as follows.
 A simple analytic continuation
argument shows that $X^0\sim \psi$ at the orbifold point 
$\psi=0$. As there is no singularity in ${\cal F}_1$ at this point,
 only the combination 
$\frac{X^0}{\psi}$ can appear in ${\cal F}_1$. 
Furthermore, we use the 
universal behavior at the conifold $$\Delta=(1-\psi^6)=0$$ and obtain, using (\ref{holFFF1}),    
\begin{equation}  
{\cal F}_1=\frac{423}{4}\log\left(\frac{X^0}{\psi}\right)+\log\left( \frac{\partial \psi}{\partial t}\right)-
{1\over 24} \log(\Delta)\ .
\end{equation} 
As a consistency check, equation (\ref{limitinfinity}) is fulfilled. A further 
consistency check is vanishing of the integer
invariants
 $$n_{1,1}=n_{1,2}=0$$ 
as expected from geometrical considerations. 
Finally, the integrality of $n_{1,d}$, which we have checked to $d=100$,
is highly non-trivial. The values for the first few  $n_{1,d}$ are listed in 
Table 2.
We also report a few of the meeting invariants as they have an interesting interpretation
as BPS bound states at threshold in Table 3.

\begin{table}
\begin{centering}
\footnotesize{
\begin{tabular}{|r|rr|}
\hline
d &g=0 &g=1\\
\hline
1&60480 & 0\\ 
2&440884080 & 0\\ 
3& 6255156277440& 2734099200\\ 
4&117715791990353760 &  387176346729900 \\ 
5& 2591176156368821985600& 26873294164654597632\\ 
6& 63022367592536650014764880&  1418722120880095142462400\\ 
7& 1642558496795158117310144372160& 65673027816957718149246220800\\ 
8&45038918271966862868230872208340160 & 2828627118403192025358734275898400 \\ 
\hline
\end{tabular}}
\caption{Integer invariants $n_{0,d}$ and $n_{1,d}$  for $X_6$ }
\end{centering} 
\end{table}

\begin{table}
\begin{centering}
\footnotesize{
\begin{tabular}{|c|rrr|}
\hline
$m_{d_1,d_2}$ &$d_2=1$ &2&3\\
\hline
$d_1=1$&15245496000&111118033656000& 1576410499948536000\\ 
2&& 809911567810170000&11490828530432030136000\\ 
3& && 163029083563567893374136000 \\ 
\hline
\end{tabular}}
\caption{Meeting invariants for $X_6$ }
\end{centering} 
\end{table}

\subsection{Quintic fibrations over $\mathbb{P}^1$}
Genus 0 Gromov-Witten invariants for multiparameter Calabi-Yau 4-folds
have been 
calculated  in \cite{Klemm:1996ts,Mayr:1996sh}. We  
determine the genus 1 Gromow-Witten invariants and test the  
integer expansion of $F_1$ for two such cases.

We first consider the quintic fibration 
over $\mathbb{P}^1$ realized as the resolution of the
degree $10$ orbifold hypersurface  
$X_{10} \subset \mathbb{P}^5(1,1,2,2,2,2)$. 
The  non-vanishing Hodge numbers are 
$$h_{0,0}=h_{4,0}=1,\ \  h_{11}=2,\ \  h_{31}=1452$$
up to symmetries of the Hodge diamond.

We introduce the divisor $F$ associated to 
the linear system  generated by monomials of degree 2. For
example, a 
representative would be $x_3=0$ yielding a degree 10 hypersurface in   
$\mathbb{P}^5(1,1,2,2,2)$. 
The 
dual curve to $F$ lies as a degree 1 curve 
in the quintic fiber with size $t_1$. Another
divisor $B$ is associated to the linear system  
generated by monomials of degree 
1. Since $B$ lies in a linear pencil of quintic fibers,
$B^2=0$. The dual curve is the base $\mathbb{P}^1$ with size $t_2$.  
  
We calculate the classical intersection data by toric geometry as follows  
\begin{equation}
\begin{array}{rl} 
F^4&=10,\ \  F^3B=5, \qquad  \ds{\int_M c_2 \wedge F^2}=\ds{110,\ \ \int_M c_2 \wedge B F=50,}\\[ 3mm]
\ds{\int_M c_4}&=\ds{2160,\quad \int_M c_3\wedge F=-200,\quad \int_M c_3\wedge B=-410 \ .} 
\end{array}
\label{topQFF}
\end{equation} 
By Batyrev's construction the mirror is given also as an degree 
10 
 hypersurface in $\mathbb{P}(1,1,2,2,2,2)/(\mathbb{Z}_{10} \mathbb{Z}_{5}^4)$,
\begin{equation} 
\sum_{k=1}^2 x_k^{2(n+1)}+\sum_{k=3}^{n+2} x_k^{(n+1)} -2 \phi \prod_{i=1}^{2} x_i^{n+1}-
\psi \prod_{k=1}^{n+2} x_i\ 
\end{equation}
with\footnote{These formulas apply to $n$-dimensional degree $2(n+1)$ hypersurfaces 
in $\mathbb{P}(1,1,2^n)$.} $n=4$.
We derive the Picard-Fuchs operators as 
\begin{equation}
\begin{array}{rl} 
{\cal L}_1&= \theta_1^{n-1}\,\left(2\theta_2-\theta_1 \right)- (n+2) 
\left[\prod_{k=1}^{n+1} ((n+2) \theta_1-k)\right] z_1\\ [ 2 mm]
{\cal L}_2&= \theta_2 - (2 \theta_2-\theta_1 -1  )(2 \theta_2-\theta_1 -2 )z_2 \ ,
\end{array}
\end{equation}
where $\theta_i=z_i \frac{{\rm d}}{{\rm d} z_i}$ with 
$z_1=\frac{\phi}{(-(n+2) \psi)^{n+2}}$ and $z_2=\frac{1}{(2 \phi)^2}$. The system 
has one conifold discriminant $\Delta_{con}$  and a `strong coupling'  
discriminant $\Delta_{s}$ at 
\begin{equation}
\Delta_{con}=1-(\psi^{n+2}- \phi)^2 ,\qquad \Delta_s=1-\phi^2\ .
\end{equation}

Let us now turn to the calculation of $n_{0,\beta}(S_i)$  and $n_{0,\beta}(c_2)$. 
We denote by $A^{(1)}_1=J_F$ and $A^{(1)}_2=J_B$ the harmonic $(1,1)$-forms dual to 
$F$ and $B$. We chose further a basis $A^{(2)}_1=J_F^2$  and 
$A^{(2)}_2=J_B J_F $ for the vertical subspace of $H^{2,2}(M)$. Toric geometry implies that 
the latter can be obtained from the leading $\theta$-polynomials $\hat {\cal L}_i$ of the Picard-Fuchs 
operators. More precisely the subspace is spanned  by the degree two elements of the 
graded multiplicative ring 
\begin{equation} 
{\cal R}=\mathbb{C}[A^{(1)}_1,\ldots,A^{(1)}_{h^{11}}]/
{\rm Id} (\hat {\cal L}_i|_{\theta_i\rightarrow A^{(1)}_i}) \ ,
\end{equation}   
where the $\hat {\cal L}_i$ are the formal limits 
$\hat {\cal L}_i=\lim_{z_i\rightarrow 0}{\cal L}_i$ of the Picard-Fuchs operators.

Following \cite{Hosono:1994ax,Klemm:1996ts,Mayr:1996sh}, we calculate the genus 0
quantum cohomology intersection  
\begin{equation} 
C^{(1)}_{ij\alpha}=\int_{M} A^{(1)}_i\wedge  A^{(1)}_j\wedge  A^{(2)}_\alpha+ {\rm instanton  
\ corrections.}
\label{cijal}
\end{equation} 
in the $B$ model as follows. Using the flat coordinates 
$$t_i=\frac{X^i}{X^0}=\frac{1}{2 \pi i}\log(z_i)+{\cal O}(z)$$ 
for the identification 
of the $B$-model structure at the point 
of maximal unipotent monodromy with the A-model structure \cite{Hosono:1994ax}, 
we find   
\begin{equation}
C^{(1)}_{ij\alpha}= \partial_{t_i} \partial_{t_j} \frac{\Pi^{(2)}_\alpha}{X^0}\ ,
\end{equation}
where given an  $A^{(2)}_\alpha$  the dual period $\Pi^{(2)}_\alpha$ is specified by the leading  
quadratic behavior in the logarithms as\footnote{Note that admixtures of periods with lower 
leading logarithmic behavior does not affect $C^{(1)}_{ij\alpha}$ 
due to the derivatives in (\ref{cijal}).} 
\begin{equation}
\frac{\Pi^{(2)}_\alpha}{X^0}=
\frac{1}{2} \sum_{ij} \int_{M} A^{(1)}_i\wedge  A^{(1)}_j\wedge  A^{(2)}_\alpha  \times  
\frac{\log(z_i) \log(z_j) }{(2 \pi i)^2}+{\cal O}(z)\ .
\end{equation}
For example, using (\ref{topQFF}),
 the period $\Pi_{c_2}$, whose expansion in $q_i$ yields the $n_{0,\beta}(c_2)$, 
is specified by the leading logarithmic behavior 
$$\frac{\Pi_{c_2}}{X^0}=\frac{1}{(2 \pi i)^2} (55 \log(z_1)^2 + 
50 \log(z_1)\log(z_2) ).$$  With this information,
 we calculate the invariants
 $$n^i_{0,\beta}= n_{0,\beta}(A^{(2)}_i)$$ 
as well as  $n_{0,\beta}(c_2)$.  

Finally, for the genus 1 Gromov-Witten invariants, we obtain
\begin{equation}  
{\cal F}_1=86 \log\left(\frac{X^0}{\psi}\right)+
\log\det\left( \frac{\partial (\psi,\phi)}{\partial (t_1,t_2)}\right)
-{1\over 24} \log(\Delta_{con})
-{7\over 24} \log(\Delta_s)
 \ .
\label{QFF}
\end{equation} 
The only difference in the calculation for the sextic is 
the  behavior at $\Delta_s$ must be determined. The latter
determination is made by imposing 
(\ref{limitinfinity}) with $\int_M c_3\wedge B=-410$. 
The second 
condition in (\ref{limitinfinity}) is a check. At $\Delta_s=0$ 
we have divisor collapsing, which is an $\mathbb{P}^1$ fibration 
over the degree 5 hypersurface in  $\mathbb{P}^3$. Integer as well as 
meeting invariants are listed in Tables 4 and 5.

\begin{table}
\begin{centering}
\footnotesize{
\begin{tabular}{|c|rrrrr|}
\hline
\!\!\!\! $n^1_{0,\beta}$\!\!\!\!  &$d_2=0$ &1 &2&3&4\\
\hline
$d_1=\! 0\!$& *& 0& 0& 0& 0 \\ 
1& 12250& 12250& 0& 0& 0 \\ 
2& 6462250& 35338750& 6462250& 0& 0 \\
3&  5718284750& 85125750000& 85125750000& 5718284750& 0\\
4& \!\!\!\!\ 6349209995000&\!\!\!\!192339896968750& \!\!\!\!\!\!507648446407500& \!\!\!\!\!\!192339896968750& \!\!\!\!\!\!6349209995000\\
\hline
\! \!\!\!$n^2_{0,\beta}$ \!\!\!\!  &$d_2=0$ &1 &2&3&4\\
\hline
$d_1\!=\!0$&  *& 5& 0& 0& 0\\
1&  2875& 9375& 0& 0& 0\\
2& 1218500& 17669375& 5243750& 0& 0\\
3&  951619125& 34150175000& 50975575000& 4766665625& 0\\ 
4&  \!\!\!\!\!\!969870120000& \!\!\!\!\!\!66623314796875& \!\!\!\!\!\!253824223203750& \!\!\!\!\!\!125716582171875& \!\!\!\!\!\!5379339875000\\
\hline
\!\!\!\! $n_{1,\beta}$ \!\!\!\!  &$d_2=0$ &1 &2&3&4\\
\hline
$d_1\!=\!0$&  *& 0& 0& 0& 0\\ 
1&  0& 0& 0& 0& 0\\
2& 0& 0& 0& 0& 0\\
3& -2768250& 7297250& 7297250& -2768250& 0\\
4& -17325370250& 90447173500& 699252105750& 90447173500& -17325370250\\
\hline\end{tabular}}
\caption{Integer invariants for the resolution of
 $X_{10}\subset \mathbb{P}^5(1,1,2,2,2)$.}
\end{centering} 
\end{table}

\begin{table}
\begin{centering}
\footnotesize{
\begin{tabular}{|c|rrrrrrr|}
\hline
$m_{\beta_1,\beta_2}$ &(0, 1)& (0, 2)& (1, 0)& (1, 1)& (1, 2)& (2, 0)& (2, 1)\\
\hline
(0, 1)&-10& -10& 6500& 0& 0& 4025250& 4025250\\ 
(0, 2)&& 10& 0& 0& 0& 0& 0\\
(1, 0)&&&10781250& 19237750& 0& 5310625000& 43309206500\\
(1, 1)&&&&10768250& 0& 10532668750& 43309206500\\
(1, 2)&&&&& 0& 0& 0\\
(2, 0)&&&&&&2555792968750& 22836787744000\\
(2, 1)&&&&&&& 124882678630250\\
\hline
\end{tabular}}
\caption{Meeting invariants for the resolution of
 $X_{10}\subset \mathbb{P}^5(1,1,2,2,2)$. }
\end{centering} 
\end{table}

It is interesting to compare the above quintic fibration 
with a different quintic fibration  given by hypersurface of  bidegree $(2,5)$, 
$$X_{2,5}\in \mathbb{P}^1\times \mathbb{P}^4\ . $$
The Hodge diamond of the second fibration is the same as the previous
case.
The divisors $B$ and $F$ correspond to the pull-backs of the
hyperplane classes
on $\mathbb{P}^1$ and $\mathbb{P}^4$ respectively. We use the same 
basis for the vertical subspace of $H^{2,2}(M)$ as before 
$A^{(2)}_1=J_F^2$  and $A^{(2)}_2=J_B J_F $.
Due to the different fibration structure, the 
topological data differ from the previous case: 
\begin{equation}
\begin{array}{rl}
\ \  F^4&=2,\quad F^3B=5, \qquad \ds{\int_M c_2 \wedge F^2}=\ds{70,\quad \int_M c_2 \wedge B F=50,}  \\[3 mm]
\ \  \ds{\int_M c_4}&=\ds{2160,\quad \int_M c_3\wedge F=-200,\quad \int_M c_3\wedge B=-330 \ .} \\ 
\end{array}
\label{topQF2F}
\end{equation} 
We derive the Picard-Fuchs equations in the standard large volume variables $z_1$ and
$z_2$~\cite{Hosono:1994ax} as 
\begin{equation} 
\begin{array}{rl} 
{\cal L}_1&= \ds{\theta_1^3\,\left(5 \theta_1- 2 \theta_2 \right)-  
5^5 z_1\prod_{k=1}^{4} (5 \theta_1+2 \theta_2 + k) + 4 z_2 (5 \theta_1+2 \theta_2 + 1)}\\ [ 1 mm]
{\cal L}_2&= \ds{\theta_2 - z_2 \prod_{k=1}^2 (5 \theta_2+2\theta_1+k  )} \ .
\end{array}
\end{equation} 
The system has only one conifold discriminant   
\begin{equation}
\Delta=(1-x_1^2)- 5x_2(1+4 x_1) + 10 x_2^2(1-x_1)   -x_2^3(10-5 x_2+x_2^2)\ .
\end{equation}
We have
introduced rescaled variables
 $x_1=5^5 z_1$ and $x_2=2^2 z_2$. Here, we know  
no further regularity conditions in the interior of the moduli space. 
Therefore, we simply 
impose (\ref{limitinfinity}) with 
$$\int_M c_3\wedge B=-330 \ \text{and} \  
\int_M c_3\wedge B=-200.$$
 That fixes the coefficients of the $\log(z_1)$ 
and $\log(z_2)$ terms in the most general ansatz of the ambiguity
\begin{equation}  
{\cal F}_1=86 \log\left(X^0\right)+
\log\det\left( \frac{\partial (z_1,z_2)}{\partial (t_1,t_2)}\right)-\frac{1}{24}\log(\Delta)
+{51\over 4} \log(z_1)
+{22\over 3} \log(z_2)\ .
\label{QF2F}
\end{equation}

The integer invariants listed in Table 6 are compatible with the previous
quintic fibration --- we get the same invariants in the fiber direction, 
as expected.

\begin{table}
\begin{centering}
\footnotesize{
\begin{tabular}{|c|rrrr|}
\hline
\!\!\!\! $n^1_{0,\beta}$\!\!\!\!  &$d_2=0$ &1 &2&3\\
\hline
$d_1=\! 0\!$&      *&      0&      0&      0\\ 
1&              9950& 171750& 609500& 609500 \\ 
2& 5487450& 533197250& 9651689750& 63917722000 \\
3& 4956989450& 1342522028500& 64483365881000& 1152361680367750\\
4& \!\!\!\! 5573313899000&\!\!\!\! 3120681190272750& \!\!\!\!\!\! 301443864603401500& \!\!\!\!\!\ 10812807897775185750 \\
\hline
\! \!\!\!$n^2_{0,\beta}$ \!\!\!\!  &$d_2=0$ &1 &2&3\\
\hline
$d_1\!=\!0$&       *&          125&              0&                  0\\
1&              2875&       195875&        1248250&            1799250\\
2&           1218500&    369229625&    10980854250&       101591346500\\
3&         951619125& 713334157250& 53873269172000&   1308427978728875\\ 
4&  \!\!\!\!\!\!969870120000& \!\!\!\!\!\! 1390949237651750& \!\!\!\!\!\! 205222409245164750& \!\!\!\!\!\! 
9819953566670512000\\
\hline
\!\!\!\! $n_{1,\beta}$ \!\!\!\!  &$d_2=0$ &1 &2&3\\
\hline
$d_1\!=\!0$&  *& 0& 0& 0\\ 
1&            0& 0& 0& 0\\
2&            0& 0& 0& 0\\
3& -2768250& 218986250& 82508848750& 2759605738750\\
4&  -17325370250& 2510820252500& 1468762788741875& 94873159058300000\\
\hline\end{tabular}}
\caption{Integer invariants for  $X_{2,5}\subset \mathbb{P}^1 \times \mathbb{P}^4$.}
\end{centering} 
\end{table}

\subsection{Elliptically fibered Calabi-Yau 4-folds}
\label{wp}
A simple 
elliptic fibration over $\mathbb{P}^3$ compactifies the local model in 
Section \ref{lcn}. Consider the resolution of the 
degree 24 orbifold hypersurface
$$X_{24}\subset \mathbb{P}^5(1,1,1,1,8,12)$$
in
weighted projective space. The genus 0 invariants 
have be calculated in \cite{Klemm:1996ts}.
 The
resolution has the following non-vanishing Hodge numbers
$$h_{0,0}=h_{4,0}=1,\ \ h_{11}=2,\ \ h_{31}=3878,\ \  
h_{22}=15564.$$
up to symmetries.

We introduce the linear system $B$ generated by linear polynomials
in the four degree 1 variables.
The linear system maps $X_{24}$ to $\mathbb{P}^3$ 
with fibers given by elliptic curves. 
%$$x_5^2+x_6^3=constant.$$ 
The second linear system $E$ is generated by polynomials of degree 4.
The curve dual to $E$ is a curve extending 
over the fiber $E$ with size denoted by $t_1$. The curve dual 
to $B$ is a degree one curve in $\mathbb{P}^3$ with size  
denoted by $t_2$.
The intersections of the divisors are
\begin{equation} 
E^4=64,\quad E^3B=16,
\quad E^2 B^2=4,\quad EB^3=1, \quad B^4=0\ .
\end{equation}
Further topological data are
\begin{equation}
\begin{array}{rl} 
\ \  \ds{\int_M c_4}&=\ds{23328,\quad 
\int_M c_3\wedge B=-960,\quad 
\int_M c_3\wedge E=-3860,} \\ [ 3mm]
\ds{\int_M c_2 \wedge B^2}&=\ds{48,\quad \ \ 
\int_M c_2 \wedge B E=182,\quad \ \ 
\int_M c_2 \wedge E^2=728.}
\end{array}
\label{topEFF}
\end{equation}

The mirror family is likewise given by an hypersurface of degree 24 in 
$\mathbb{P}(1,1,1,1,8,12)/(\mathbb{Z}_{24}^3)$
\begin{equation} 
\sum_{k=1}^n x_1^{6n}+x_{n+1}^2+x_{n+2}^3-n\phi \prod_{k=1}^n x_i^{2 n}-6 n \psi \prod_{i=1}^{n+2} x_i\ 
\end{equation}
with\footnote{These formulas apply to $n$-dimensional degree $6n$ hypersurfaces 
in $\mathbb{P}(1^n,2 n,3n)$.} $n=4$.
We derive the Picard-Fuchs operators as 
\begin{equation}
\begin{array}{rl} 
{\cal L}_1&= \theta_1(\theta_1-n \theta_2)-12 (6\, \theta_1-5)(6\, \theta_1-1) z_1\\ [ 2 mm]
{\cal L}_2&= \theta_2^n -\prod_{k=1}^n (n\, \theta_2-\theta_1-k) z_2\ ,
\end{array}
\end{equation}
where $\theta_i=z_i \frac{{\rm d}}{{\rm d} z_i}$ with 
$z_1=\frac{n \phi}{(n \psi)^6}$ and $z_2=\frac{(-1)^n}{(n \phi)^n}$.
The system has two conifold discriminants
\begin{equation}
\Delta_1=1-\phi^n,\qquad \Delta_2=1-\tilde\phi^n\ ,
\end{equation}
where we defined $\tilde \phi=\psi^6-\phi$. The solutions to the Picard-Fuchs equations
can be obtained similarly as in Section
 \ref{sextic} using the methods outlined in \cite{Klemm:1996ts}.
For example the holomorphic solution at the point of maximal unipotent monodromy is given by
\begin{equation}
X^0=\sum_{k_1=0,k_2=0}^\infty \frac{ (6k_1)! (n k_2)!}{  (2 k_1)! (3 k_1)! k_1! (k_2!)^n} z_1^{k_1} z_2^{k_2} \ .
\end{equation}

The considerations, which lead to the expression of $F_1$ in the holomorphic limit, are very
similar to those of  Section \ref{sextic}, 
\begin{equation}  
{\cal F}_1=928 \log\left(\frac{X^0}{\psi}\right)+
\log\det\left( \frac{\partial (\psi,\phi)}{\partial (t_1,t_2)}\right)+3\log(\psi)-
{1\over 24} \sum_{i=1}^2\log(\Delta_i) \ .
\label{FHEF}
\end{equation} 
A new feature  here is  $$\lim_{\psi\rightarrow 0}
\log\det\left( \frac{\partial (\psi,\phi)}{\partial (t_1,t_2)}\right)\sim\psi^{-3},$$
 as  is 
shown by simple analytic continuation of $X^0$ and the two logarithmic solutions to $\psi=0$. 
To maintain the expected regularity at $\psi=0$, we have to add the explicit 
$3\log(\psi)$ term to the holomorphic ambiguity. As a check of the result (\ref{FHEF}),
 we 
note again that (\ref{limitinfinity}) with (\ref{topEFF}) is fulfilled.

We chose further a basis $A^{(2)}_1=\frac{1}{17}(4 J_E^2+J_E J_B)$  and 
$A^{(2)}_2=J_B^2$ and 
 calculate as before the genus 0 and genus 1 invariants. 
As a consistency 
check we note that scaling the size of the elliptic fiber $t_1$ to infinity 
leaves us precisely with the 
${\cal O}(-4)\rightarrow \mathbb{P}^3$ geometry. 
The corresponding invariants are listed in Table 7 .

\begin{table}
\begin{centering}
\footnotesize{
\begin{tabular}{|c|rrrrr|}
\hline
\!\!\!\! $n^1_{0,\beta}$\!\!\!\!  &$d_2=0$ &1 &2&3&4\\
\hline
\!$d_1=0$\!&  0 &    0&      0&        0& 0\\ 
1& 960& 5760& 181440& 13791360& 1458000000\\ 
2& 1920& -1817280& -98640000& -10715760000& -1476352644480\\
3& 2880& 421685760& 29972448000& 4447212981120& 783432258136320\\
4& 3840& \!\!\!\!\!\! 2555202430080 &\!\!\!\!\!\! -6353500619520& \!\!\!\!-1273702762398720& \!\!\!\!\!\!\!-285239128072550400\\
\hline
\! \!\!\!$n^2_{0,\beta}$ \!\!\!\!  &$d_2=0$ &1 &2&3&4\\
\hline
\!$d_1=0$\!&  0& -20& -820& -68060& -7486440\\ 
1&  0& 7680& 491520& 56256000& 7943424000\\ 
2&  0& -1800000& -159801600& -24602371200& -4394584496640\\
3&  0& 278394880& 35703398400& 7380433205760& 1662353371955200\\
4&  0& \!\!\!  623056099920 &\!\!\!\!\! -6039828417600& \!\!\!\!\! -1683081588149760& \!\!\!\!\!\!-478655396625235200\\
\hline
\!\!\!\! $n_{1,\beta}$ \!\!\!\!  & $d_2=0$&1 &2&3&4\\
\hline
\!$d_1=0$\!&   0& 0& 0& 11200& 3747900\\ 
1& -20& -120& -3780& -7852120& -3536410200\\ 
2& 0& 45720& 2245680& 2858334000& 1724679193440  \\
3& 0& -10662240& -719326800& -719497580160& -573686979645680 \\
4& 0& 1638152760& 160844654520& 140278855296640&\!\!\!\!\!\! 145314212874711600\\
\hline
\end{tabular}}
\caption{Integer invariants for the resolution of $X_{24}$.}
\end{centering} 
\end{table}

%\pagebreak

\vspace{+10 pt}
\noindent
Department of Physics \\
Univ. of Wisconsin \\
Madison, WI 53706, USA\\
aklemm@physics.wisc.edu \\

\vspace{+10 pt}
\noindent
Department of Mathematics\\
Princeton University\\
Princeton, NJ 08544, USA\\
rahulp@math.princeton.edu

\end{document}